\documentclass[11pt]{amsart}     
\usepackage{a4wide,graphicx}
\usepackage{color}
\usepackage{caption}
\usepackage{comment}
\usepackage{enumerate}
\usepackage{upgreek}
\usepackage[left=1in, right=1in]{geometry}
\usepackage{bbm,yfonts}
\usepackage{hyperref,amssymb} 

\let\pa\partial  
  
\let\eps\varepsilon

\newcommand{\R}{{\mathbb R}} 

\newcommand{\diver}{\operatorname{div}}
\newcommand{\even}{\mathrm{e}}
\newcommand{\odd}{\mathrm{o}}

\newcommand{\T}{\mathsf T}

\newcommand{\Bb}{\operatorname{\mathbb{B}}}
\newcommand{\Db}{\operatorname{\mathbb{D}}}
\newcommand{\dd}{{\mathrm{d}}}

\newcommand{\sym}{\operatorname{sym}}

\newcommand{\V}{\mathcal{V}}
\newcommand{\E}{\mathcal{E}}
\newcommand{\Diss}{\mathcal{D}}
\newcommand{\W}{\mathcal{W}}

\DeclareMathOperator*{\esssup}{ess\,sup}
\newcommand{\vphi}{\boldsymbol\varphi}

\newcommand{\vro}{\varrho}
\newcommand{\bv}{\boldsymbol v}
\newcommand{\bu}{\boldsymbol u}
\newcommand{\bs}{\boldsymbol}
\newcommand{\bx}{\boldsymbol x}
\newcommand{\by}{\boldsymbol y}
\newcommand{\bz}{\boldsymbol z}
\newcommand{\va}{\mathsf v}
\newcommand{\ua}{\mathsf u}
\newcommand{\paa}{\mathsf p}
\newcommand{\res}{\mathsf r}

\newcommand{\material}{\Omega}
\newcommand{\Hper}{H_{\mathrm{per}}}
\newcommand{\Rey}{\operatorname{Re}}
\newcommand{\el}{\operatorname{\mathbb{L}}}
\newcommand{\Dom}{\mathrm{Dom}}

\newtheorem{theorem}{Theorem}[section]

\theoremstyle{definition}
\newtheorem{definition}{Definition}[section]

\theoremstyle{remark}
\newtheorem{remark}{Remark}[section]

\begin{document}

\title[A review on reduced FSI problems]{A review on rigorous derivation of reduced models\\ 
for fluid - structure interaction systems}
 
\author{Mario Bukal}
\author{Boris Muha}
\address{University of Zagreb,
Faculty of Electrical Engineering and Computing\newline
Unska 3, 10000 Zagreb, Croatia}
\email{mario.bukal@fer.hr}
\address{University of Zagreb,
Faculty of Science, Department of Mathematics,
Bijeni\v cka cesta 30, 10000 Zagreb, Croatia}
\email{borism@math.hr}

\thanks{This work has been supported in part by the Croatian Science
Foundation under projects 7249 (MANDphy) and 3706 (FSIApp).}

\keywords{thin viscous fluids, fluid - structure interaction problems, reduced models,
sixth-order thin-film equations, error estimates.}

\subjclass[2010]{35M30, 35Q30, 35Q74, 76D05, 76D08}
%
%
%
\begin{abstract}
In this paper we review and systematize the mathematical theory on justification of sixth-order 
thin-film equations as reduced models for various fluid - structure interaction systems in which
fluids are lubricating underneath elastic structures. Justification is based on careful 
examination of energy estimates, weak convergence results of solutions of the original 
fluid - structure interaction systems to the solution of the sixth-order 
thin-film equation, and quantitative error estimates which provide even
strong convergence results.
\end{abstract}

\date{\today}   
\maketitle

%


\section{Introduction}

In 1886 Reynolds derived the fundamental equation of the lubrication approximation \cite{Rey86},
which serves until nowadays in many engineering applications. 
It is an elliptic equation for the pressure distribution in relatively thin viscous fluid
in laminar flow between two rigid parallel plates in a relative motion of constant velocity.
The equation can be understood as a reduced model for the basic Navier-Stokes equations 
describing the fluid motion. 
In Section \ref{sec:heuristic} we outline the heuristic derivation of the Reynolds equation.
However, often in nature and in engineering applications, those ``parallel plates'' are not rigid,
but have their own dynamics which also affects the fluid quantities. Such systems are widely known
as fluid - structure interaction (FSI) systems. They appear in medicine, in particular in modelling 
of cardio-vascular systems \cite{BGN14,BCMG16}, then in aero-elasticity \cite{Dow15}, 
marine engineering \cite{YGJ17}, etc.~and as interesting systems of partial differential
equations, they also gained a huge attention in 
applied mathematics comunity.  

In this review paper we are particularly interested in FSI systems in which fluids are 
lubricating underneath elastic structures. Such models describe for instance the growth of magma 
intrusions \cite{LPN13, Mic11}, the fluid-driven opening of fractures in the Earth's crust 
\cite{BuDe05, HBB13}, subglacial floods \cite{DJBH08,TsRi12},
the passage of air flow in the lungs \cite{HHS08}, and the operation of vocal cords \cite{Tit94}.
They are also inevitable in engineering, for example 
in manufacturing of silicon wafers \cite{HuSo02, King89}, 
suppression of viscous fingering \cite{PIHJ12, PJH14} and  
in an emerging area of microfluidics \cite{LBS05,HoMa04,TaVe12}
with particular applications to so called lab-on-a-chip technologies \cite{SSA04, DaFin06}.
Describing such systems with their true physical models: Navier-Stokes equations coupled with
elasticity equations both on relatively thin domains, is inappropriate from analytical
and numerical points of view, and thus, inappropriate for engineering applications.
Therefore, reduced and simplified models are sought which will maintain the essence of the
fluid - structure interaction. Depending on the original problem at hand, many such reduced
models have been derived, especially in engineering literature. We emphasize at this point
that our aim is not to cover all those examples, that would be impossible, but we consider 
several important model examples and concentrate on rigorous derivation of the reduced models 
and identification of the necessary scaling assumption in system paramaters, which ``sees'' 
the interaction between the two subsystems.   

Using asymptotic expansion techniques several reduced models of Biot-type describing the flow through a 
long elastic axially symmetric channel have been derived in \cite{CanMik03,MikGuCan07,TamCanMik05}.
A rigorous justification of those reduced model by means of weak convergence results and the corresponding 
error estimates was provided  in \cite{CanMik03}. Periodic flow in thin channel with visco-elastic walls 
was analyzed in \cite{PaSt06}, where starting from a linear 2D (fluid)/1D (structure) FSI model and 
under particular assumption on the ratio of the channel 
height and rigidity of the wall, a linear sixth-order thin-film equation describing the wall 
displacement emanated as a reduced model. A similar problem has been also addressed in \cite{CuMP18}, 
and the reduced model in terms of a linear sixth-order equation arose again.
The reduced models in both papers have been justified by the corresponding weak convergence results. 

In this review paper we mainly focus on linear FSI problems, meaning that equations for both subsystems are linear
and moreover, the interaction between the fluid and the structure is realized through a fixed interface. 
The main reason for such simple FSI models is that they provide global in time existence of weak 
solutions which possess sufficient regularity for passing from the FSI to 
the reduced model in a rigorous way by means of weak convergence results. We discuss two types of the 
linear FSI problems.
First, a 3D/2D FSI problem, analogous to the FSI problem in \cite{PaSt06}, is discussed in 
Section \ref{sec:e+0}. In this model the structure is originally described by a lower-dimensional 
elasticity model. Hence, the dimension reduction only applies to the
fluid part. Additional horizontal dimension in our case does not bring any conceptual novelty.  
In contrast to the asymptotic expansion 
techniques employed in \cite{PaSt06}, our approach relies on careful examination of energy estimates. 
Based on these, quantitative a priori estimates in terms of the fluid thickness have been derived 
for weak solutions, and weak convergence results have been established. 
Moreover, these quantitative a priori estimates provide the necessary scaling assumption on model parameters 
which ensures the nontrivial reduced model. 
The same ideas have been recently employed in \cite{BuMu19}, where the authors analyzed a linear 
3D/3D FSI problem, in which the simultaneous dimension reduction in the fluid and the structure has been 
performed and again a linear sixth-order thin-film equation has been derived as the reduced model 
under particular scaling assumptions. 
We briefly report on this problem and results in Section \ref{sec:e+h}.
Finally, we address the physically most relevant nonlinear FSI problem in which equations of fluid 
motion are nonlinear and the fluid domain is also unknown in the system. Unfortunately, the 
wellposedness results for such problems are very scarce, hence, for the beginning we restrict our
analysis to a 2D/1D FSI problem for which the existence of global in time strong solutions is 
available from \cite{grandmont2016existence}. Rigorous derivation of the nonlinear sixth-order 
thin-film equation, which is known in the engineering literature \cite{HoMa04, HBB13, LPN13}, 
as the reduced model for this FSI problem is still work in progress
\cite{BuMu20}, and only main ideas are outlined in Section \ref{sec:nonlin}.

\section{Heuristic derivation of reduced models }\label{sec:heuristic}
We start with the heuristic derivation of reduced models, which is essentially based on the lubrication approximation 
in the fluid part and the assumption of the pressure balance on the interface.

Let us consider the Navier-Stokes system in a time-dependent domain
$\Omega_\eta = \{(x,z)\in (0,L)\times(0,\eta(x,t))\}\subset\R^2$, where $L>0$ is a fixed length and
$\eta(x,t)\geq\eta_0>0$ is a given uniformly positive, $L$-periodic smooth function 
which depicts the evolution of the upper boundary of $\Omega_\eta$. The fluid velocity
$\bv=(v_1,v_2)$ and pressure $p$ then satisfy
\begin{align}\label{intro:eq:NS}
\vro_f\pa_t\bv + \vro_f(\bv\cdot\nabla)\bv - \mu\Delta\bv = -\nabla p\,,
&\quad\text{on }\Omega_\eta\times(0,\infty)\,,\\
\diver \bv =0\,,&\quad\text{on }\Omega_\eta\times(0,\infty)\,.
\label{intro:eq:divfree}
\end{align}
On the bottom boundary we consider fluid moving with a constant velocity in
horizontal di\-rec\-ti\-on, i.e.~$\bv(\cdot,0,\cdot)=(v_D,0)$, and
on the top boundary the fluid moves with the velocity of the surface, i.e.~{$v_2 = \pa_t\eta$}. 
Concerning the lateral boundaries, for the completeness of the problem we may take the periodic 
boundary conditions, but for the purpose of this exposition, they are not so important.
Imposing initial velocity $\bv(\cdot,\cdot,0)=\bv_0$ completes the problem 
(\ref{intro:eq:NS})-(\ref{intro:eq:divfree}).

In many physical situations like spilled water on floor, fresh paint on a wall, or industrially
more relevant microfluidics devices and lab-on-a-chip technologies \cite{SSA04, DaFin06, LLPT15}, 
one typically has the length in 
one spatial direction (typically vertical) much smaller than the other one, i.e.~$H:=\sup_{x,t}\eta(x,t)\ll L$. 
Let us denote $\eps:=H/L$ and write down equations (\ref{intro:eq:NS}) and (\ref{intro:eq:divfree}) 
in a non-dimensional form in terms of the small parameter $\eps$. 
For that purpose we introduce new non-dimensional variables:
\begin{align*}
\hat{x} = \frac{x}{L}\,,\quad \hat{z} = \frac{z}{H}\,,\quad \hat{v}_\alpha = 
\frac{v_\alpha}{V_\alpha}\,,\ \alpha = 1,2\,, \quad \hat{p} = \frac{p}{P}\,,
\quad \hat{t}=\frac{t}{\T}\,,\quad \hat{\eta} = \frac{\eta}{H}\,,
\end{align*}
where $V_\alpha$, $P$ and $\T$ denote nominal values of fluid velocities, pressure and the 
time scale, respectively. The following calculations are standard in many fluid mechanics textbooks 
or lecture notes, see for instance \cite{Sze12}. 
Performing the above change of variables in (\ref{intro:eq:NS}) 
and (\ref{intro:eq:divfree}), and neglecting the hat notation in new variables, we find:
\begin{align*}
\frac{\vro_fV_1L}{PT}\pa_t v_1 + \frac{\vro_fV_1^2}{L}\left(v_1\pa_x v_1 
+ \frac{V_2}{V_1}\eps^{-1} v_2\pa_zv_1\right)
- \frac{\mu V_1}{LP}\left(\pa_{xx}v_1 + \eps^{-2}\pa_{zz}v_1 \right) + \pa_x p = 0\,,\\
\frac{\vro_fV_2L}{PT}\pa_t v_2 + \frac{\vro_fV_2^2}{L}\left(\frac{V_1}{V_2}v_1\pa_x v_2 
+ \eps^{-1} v_2\pa_zv_2\right)
- \frac{\mu V_2}{LP}\left(\pa_{xx}v_2 + \eps^{-2}\pa_{zz}v_2 \right) + \eps^{-1}\pa_z p = 0\,,\\
\frac{V_1}{V_2}\pa_x v_1 + \eps^{-1}\pa_z v_2 = 0\,.
\end{align*}
Notice that $\mu V_1/LP$ is a dimensionless quantity. It is customary in the lubrication
approximation regime to assume that $\mu V_1/LP = \eps^2$ and $V_2 = \eps V_1$ \cite{Sze12}.
Then the above system becomes
\begin{align}\label{intro:eq:NS-res}
\eps^2\Rey\left(\pa_t v_1 + v_1\pa_x v_1 + v_2\pa_zv_1\right)
- \eps^2\pa_{xx}v_1 - \pa_{zz}v_1 + \pa_x p = 0\,,\\
\eps^4\Rey\left(\pa_t v_2 + v_2\pa_x v_1 + v_2\pa_zv_2\right)
- \eps^4\pa_{xx}v_2 - \eps^{2}\pa_{zz}v_2 + \pa_z p = 0\,,\\
\pa_x v_1 + \pa_z v_2 = 0\,,\label{intro:eq:divfree-res}
\end{align} 
where $\Rey = \vro_fL^2/(\mu \T)$ denotes the Reynolds number. Under
assumption of $\eps\ll1$ and $\eps^2\Rey\ll 1$, and thus neglecting terms of those orders in 
(\ref{intro:eq:NS-res})-(\ref{intro:eq:divfree-res}), we arrive to the system
\begin{align}
- \pa_{zz}v_1 + \pa_x p = 0\,,\label{intro:eq:poisson}\\
\pa_z p = 0\,,\label{intro:pressy}\\
\pa_x v_1 + \pa_z v_2 = 0\,. \label{intro:divfreer}
\end{align}

Integrating equation (\ref{intro:divfreer}) with respect to $z$ and employing 
the corresponding boundary conditions we find
\begin{equation}\label{eq:v2}
\partial_t\eta  
= - \pa_x \int_0^{\eta(x,t)} v_1(x,z)\,\dd z\,.
\end{equation}
Since equation (\ref{intro:pressy}) implies that the pressure is constant in the vertical direction,
equation (\ref{intro:eq:poisson}) is a Poisson equation with respect to $z$ with Dirichlet
boundary condition. Its explicit solution is then given by,
\begin{equation}
v_1(x,z,t) = \frac12z(z-\eta(x,t))\pa_x p(x,t) + \left(1 - \frac{z}{\eta} \right)v_D\,.
\end{equation}
Integrating the last equation with respect to $z$ and utilizing the obtained expression in (\ref{eq:v2}), 
we arrive to the elliptic equation for the pressure $p$:
\begin{equation}\label{eq:reynolds}
- \pa_x\left(\frac{\eta^3}{12}\pa_x p \right) = -\pa_t\eta - \frac12\pa_x(\eta v_D)\,.
\end{equation}

Assuming that $\eta$ is a stationary profile and horizontal velocity $v_0$ is constant,
equation (\ref{eq:reynolds}) turns into the original Reynolds equation
\begin{equation}
-\pa_x\left(\eta^3\pa_x p \right) = -6v_D\pa_x\eta\,.
\end{equation}
This is a fundamental equation of elastohydrodynamics, which has been derived by Reynolds himself in
\cite{Rey86}, and much later rigorously justified in \cite{BaCh86,Cim83}.

However, in most applications the upper boundary $\eta$ is not a priori known, but it is coupled with 
pressure $p$ through another equation. 
Physically, such equation describes the balance of forces on the interface between 
the two phases and in the sequel we discuss three most common physical situations:
\begin{enumerate}[(i)]
  \item When the pressure is balanced by the gravity, then the pressure is proportional to 
  		$\eta$ \cite{HBB13}. 
  		After appropriate time rescaling, equation (\ref{eq:reynolds}) then 
  		becomes a porous medium type equation
  		\begin{equation}\label{2:pm}
			\pa_t\eta = \pa_x^2\left(\eta^4 \right) - 6\pa_x(\eta v_D)\,.
		\end{equation}
		The mathematical theory of the porous medium type equations like (\ref{2:pm}) 
		is very well developed and comprehended in the monograph \cite{Vaz07}.
  \item In the presence of dominant surface tension force, the pressure is balanced by the linearized 
  curvature, i.e.~$p \sim - \pa_x^2\eta$ \cite{Mye98}. Equation (\ref{eq:reynolds})
  		then becomes the well known thin-film equation \cite{Ber98, Mye98}
  		\begin{equation}\label{2:tf}
			\pa_t\eta = -\pa_x\left(\eta^3\pa_x^3 \eta \right) - 6\pa_x(\eta v_D)\,.
		\end{equation}
		Fourth-order thin-film equations also gained huge attention in the applied mathematics and
		engineering community. We refer to \cite{BeFr90, CDGKSZ93, ODB97} and references therein.
		Although they do not share some fundamental concepts of the second-order equations, 
		like the maximum principle, which makes the analysis of fourth-order equations more difficult, 
		they have rich mathematical structure which has been explored in numerous papers 
		(cf.~for instance \cite{GiaOtt01, DGG98, BeGr}). 
		 
  \item When the fluid phase is covered by an elastic plate in dominantly bending regime, then the 
  pressure satisfies $p \sim \pa_x^4\eta$ and equation (\ref{eq:reynolds}) reads
  		\begin{equation}\label{2:stf}
			\pa_t\eta = \pa_x\left(\eta^3\pa_x^5 \eta \right) - 6\pa_x(\eta v_D)\,.
		\end{equation}
		Contrary to the fourth-order equations, the sixth-order thin-film equations 
		which are physically as relevant as (\ref{2:tf}), did not gain a comparable attention in the literature.
		Equations of type (\ref{2:stf}) have been derived for instance in \cite{HoMa04,LPN13}.
\end{enumerate}
Besides the gravity, other potential forces like capillarity, heating, Van der Waals forces, etc.~with 
potential $\Phi(\eta)$ can be included into physical models (i)-(iii) leading to a general equation of type
\begin{equation*}
			\pa_t\eta = (-1)^{(\alpha-1)/2}\pa_x\left(\eta^3\pa_x^\alpha \eta \right) 
			+ \pa_x\left(\eta^3\pa_x\Phi'(\eta)\right)  - 6\pa_x(\eta v_D)\,
\end{equation*}
with $\alpha=1,3,5$. A review on the derivation of plethora of such fourth-order models can be found 
in \cite{ODB97}.

\section{Reduced models for fluid - structure interaction problems: a rigorous approach}  
In this section we focus on FSI problems, formally the case (iii) of the previous section. 
We will consider fluid - structure interaction problems where both phases, fluid and structure, 
are relatively thin.
Starting from an FSI problem, our aim is to rigorously justify the reduced model in terms of a
sixth-order evolution equation of type (\ref{2:stf}). 
This means to prove that solutions of the original FSI problem
converge in some sense to the solution of the sixth-order thin-film equation, and vice-versa, by 
solving the the sixth-order thin-film equation, one is able to construct an approximate solution
to the original problem.
Rigorous justification of reduced models is so far available only for linear problems due to the
global in time wellposedness for the weak solutions, and such will 
be discussed here. Linear in this context means that equations of motion for both fluid and structure 
in FSI problem are linear and moreover, the fluid domain is fixed and therefore the coupling is 
linear and is realized on the fixed fluid - structure interface. 

As a model problem we consider a three-dimensional channel of relative
height $\eps > 0$ which is filled with incompressible viscous fluid and the channel
is covered by an elastic plate of relative height $h>0$. We work in physical 3D space,
although some results in the literature are available only in 2D, but we will emphasize when 
it comes to that point. 
Assume that the problem is properly nondimensionalized and 
denote by $\material = \Omega_\eps \cup S_h\subset\R^3$ {\em the material domain}, where 
$\Omega_\eps = (0,1)^{2}\times(-\eps,0)$ denotes the fluid domain, 
and $S_h$ denotes the structure domain, which depends on the structure model at hand. 
If we describe the structure dynamics by a lower-dimensional model, for instance linear 
visco-elastic plate, then $S_h = (0,1)^{2}\equiv\omega$ and we denote such problem
as $\eps+0$ {\em problem} (cf.~Section \ref{sec:e+0}). On the other hand, if the
structure dynamics is fully described by linear elasticity equations, then $S_h = \omega\times (0,h)$
and the problem is denoted by $\eps+h$ {\em problem} (cf.~Section \ref{sec:e+h}).
The linear FSI problem is in general described by the system of partial differential equations:
\begin{align}
\vro_f\pa_t\bv - \diver\sigma_f(\bv,p) &= \bs f\,,\quad \Omega_\eps\times(0, \infty)\,,\label{1.eq:St}\\
\diver \bv &= 0\,,\quad \Omega_\eps\times(0, \infty)\,,\label{1.eq:DF}\\
\vro_s\pa_{tt}\bu + \el(\bu,\pa_t\bu) &= \bs g\,, \quad S_h\times(0, \infty)\,,\label{1.eq:EL}
\end{align}
where equations (\ref{1.eq:St}) and (\ref{1.eq:DF}) denote the Stokes system for the 
fluid velocity $\bv$ and the pressure $p$. The fluid Cauchy stress tensor is given by
$\sigma_f(\bv,p) = 2\nu \sym \nabla \bv - pI_3$, where $\nu$, $\varrho_f>0$ denote the fluid viscosity 
and density, respectively, $\sym(\cdot)$ denotes symmetric part of the matrix and $\bs f$ denotes 
the fluid external force. The structure displacement $\bu$ is described by elasticity 
equation(s) (\ref{1.eq:EL}), where
$\el$ denotes the model dependent elasticity operator, which 
can be decomposed as a sum $\el(\bu,\pa_t\bu) = \Bb\bu + \Db\pa_t\bu$ of symmetric and positive 
operators $\Bb$ and $\Db$ which correspond to elastic and viscoelastic energy of the structure, 
respectively, and are to be specified below. Coefficient $\vro_s$ is the structure density
and $\bs g$ denotes the force of fluid acting on the structure. The structure volume forces like 
for instance gravity are for simplicity excluded from our analysis.

The two subsystems need to be coupled through the
interface conditions which we literary describe as: {\em continuity of velocities} (kinematic condition)
and {\em balance of forces} (dynamic condition). Depending 
on the model at hand they will be specified below.
For simplicity of exposition we assume periodic boundary conditions 
in horizontal variables for all unknowns. 
On the bottom of the channel
we assume the standard no-slip boundary condition for the fluid velocity, $\bv = 0$, 
and the plate is free on the top boundary. The system is supplemented by 
trivial initial conditions:
\begin{align}\label{def:IC}
\bv(0)=0,\; \bu(0)=0,\; \partial_t\bu(0)=0\,,
\end{align}
although, all obtained results will also hold for nontrivial initial conditions under some additional
smallness assumptions (cf.~for instance \cite[Appendix]{BuMu19}). 
A non\-tri\-vial volume force on the structure could also be involved, again under certain 
scaling assump\-tions (cf.~again \cite{BuMu19}), but the trivial one is in fact 
motivated by applications in microfluidics \cite{SSA04}.
The previously settled framework also incorporates physically more
relevant problem which involves prescribed pressure drop between inlet and outlet of the channel,
instead of the periodic boundary conditions.
As described in \cite{PaSt06}, this is a matter of the right choice of the fluid volume force $\bs f$.
The simplified linear FSI problem (\ref{1.eq:St})-(\ref{1.eq:EL}) can be seen as a
linearization of a truly nonlinear dynamics under the assumption of small displacements \cite{YSW16}.

Let us now proceed with a formal analysis.
Testing equations (\ref{1.eq:St}) and (\ref{1.eq:EL}) with assumed smooth solutions 
$\bv$ and $\pa_t\bu$,
respectively, integrating by parts and utilizing the divergence free condition we have: for every 
$t>0$
\begin{align}
\frac{\vro_f}{2}\int_{\Omega_\eps}|\bv(t)|^2 \dd \bx
+ 2\nu\int_0^t\!\!\int_{\Omega_\eps}|\sym\nabla\bv(s)|^2 \dd\bx \dd s 
+ \frac{\vro_s}{2}\int_{S_h}|\pa_t\bu(t)|^2\dd V \nonumber\\
+ \frac12\int_{S_h}\Bb\bu(t)\cdot \bu(t)\dd V 
+ \int_0^t\!\!\int_{S_h}\Db\pa_t\bu(s)\cdot\pa_t \bu(s)\dd V\dd s
 = \int_0^t\!\!\int_{\Omega_\eps}\bs f\cdot\bv\, \dd\bx\dd s\,,\label{EE}
\end{align}
where $\dd V$ denotes the volume measure on the structure domain. Assuming that the fluid volume 
force satisfies $\|\bs f\|_{L^\infty(0,\infty;L^\infty(\Omega_\eps))} \leq C$ and employing then the 
Poincar\'e and Korn inequalities on thin domains (cf.~\cite[Proposition A.2]{BuMu19}), the right-hand
side can be estimated and leads to the basic energy estimate: for every 
$t>0$
\begin{align}
\E_k(\bv(t))
&+ \nu\int_0^t\Diss_f(\bv(s)) \dd s 
+ \E_k(\pa_t\bu(t)) 
+\E_{el}(\bu(t)) + \int_0^t\Diss_s(\pa_t\bu(s)) \dd s  \leq Ct\eps\,, \label{eq.BEI}
\end{align}
where $C>0$ from now on denotes a generic positive constant independent of $\eps$ and $t$. 
Quantities $\E_k(\bv)$ 
and $\E_k(\pa_t\bu)$ denote kinetic energy of the fluid and the structure, 
$\Diss_f(\bv)$ and $\Diss_s(\pa_t\bu)$ denote the rate of the energy dissipation of the fluid
and the structure, respectively, while $\E_{el}(\bu)$ denotes the elastic energy of the structure. 
All quantities are all easily read off from (\ref{EE}). 

Next, we briefly describe the concept of weak solutions on an abstract level, while
details regarding the specific model are addressed in respective subsections below.
The choice of appro\-priate solution 
spaces is motivated by the above energy estimate.
For the fluid velocity, this appears to be
\begin{align*}
\quad \V_F(0,T;\Omega_\eps) = L^\infty(0,T;L^2(\Omega_\eps;\R^3))\cap L^2(0,T;V_F(\Omega_\eps))\,,
\end{align*}
where $V_F(\Omega_\eps) = \left\{\bv\in H^1(\Omega_\eps;\R^d)\ 
:\ \diver\bv = 0\,,\ \bv|_{x_3=-\eps}=0\,,\ \bv \text{ is }\omega\text{-periodic} \right\}$, and
$T>0$ is a given time horizon.
Similarly, the structure function space will be
\begin{align}\label{def:str_space}
\V_S(0,T;S_h) = W^{1,\infty}(0,T;L^2(S_h;\R^3))\cap L^\infty(0,T;\Dom(\Bb^{1/2})) 
\cap H^1(0,T;\Dom(\Db^{1/2}))\,,
\end{align}
where $\Dom(\Bb^{1/2})$ and $\Dom(\Db^{1/2})$ denote domains of respective operators. 
Finally, the solution space of the coupled problem (\ref{1.eq:St})-(\ref{def:IC}) will be compound 
of previous spaces involving the kinematic interface condition (k.\,c.) as a constraint:
\begin{align}\label{def:sol_space} 
\V(0,T;\Omega) = \big\{(\bv,\bu)\in \V_F(0,T;\Omega_\eps)&\times \V_S(0,T;S_h)\ : 
\text{ k.\,c.~holds for a.e. }t\in(0,T)\big\}\,.
\end{align}
Now we can state the definition of weak solutions to our problem in the sense of Leray and Hopf.
\begin{definition}\label{def:weak_sol}
We say that a pair $(\bv^\eps,\bu^h)\in \V(0,T;\Omega)$ is a \emph{weak solution}
to the linear FSI problem (\ref{1.eq:St})--\eqref{def:IC}, 
if the following variational equation holds in $\mathcal{D}'(0,T)$:
\begin{align}
&\vro_f\frac{\dd}{\dd t}\int_{\Omega_\eps}\bv^\eps\cdot\bs\phi\,\dd \bx
-\vro_f\int_{\Omega_\eps}\bv^\eps\cdot\pa_t\bs\phi\,\dd \bx \nonumber
  + 2\nu\int_{\Omega_\eps}\sym\nabla\bv^\eps:\sym\nabla\bs\phi\,\dd \bx 
  + \vro_s\frac{\dd}{\dd t}\int_{S_h}\bu^h\cdot\pa_t\bs\psi\,\dd S \\
& - \vro_s\int_{S_h}\pa_t\bu^h\cdot\pa_t\bs\psi\,\dd S + \label{eq:weak_sol}
\int_{S_h}\el(\bu^h,\partial_t\bu^h)\cdot\bs\psi\,\dd S
= \int_{\Omega_\eps}\bs f\cdot\bs\phi \,\dd \bx
+ \int_{S_h}\bs g\cdot\bs\psi \,\dd S
\end{align}
for all $(\bs\phi,\bs\psi)\in \W(0,T;\Omega)$, where
\begin{align*}
\W(0,T;\Omega) = \big\{ (\bs\phi,\bs\psi)\in 
C^1\left([0,T]; V_F(\Omega_\eps)\times V_S(S_h)\right)\, :\,
 \left.\bs\phi(t)\right|_{\omega} = \left.\bs\psi(t)\right|_{\omega} 
 \ \forall\, t\in [0,T] \big\}\,
\end{align*}
denotes the space of test functions and $V_S(S_h) = \Dom(\Bb^{1/2})\cap \Dom(\Db^{1/2})$.
Moreover, $(\bv^\eps,\bu^h)$ verifies the energy 
inequality (\ref{eq.BEI}).
\end{definition}

For the existence of weak solutions one typically employs the Galerkin method and formal
estimate (\ref{eq.BEI}) provides crucial a priori estimates 
needed for the construction of a unique weak solution. 
The pressure in the system is treated in a standard manner, but unlike in the Stokes system solely, 
where the pressure is determined up to a function of time, 
here in the case of the full FSI problem the pressure is unique. 
This is a consequence of the fact that in the Stokes system the boundary (wall) is assumed to 
be rigid and therefore cannot ``feel'' the pressure, while in the present case elastic wall 
feels the pressure. 

In the subsequent sections we address FSI problems (\ref{1.eq:St})-(\ref{def:IC}) 
depending of the choice of the structure model, i.e.~chioce of the elasticity operator $\el$.

\subsection{Linear $\eps + 0$ problem}\label{sec:e+0}
First we discuss the FSI problem in which the elastic plate covering the fluid channel 
is already treated as a lower-dimensional object.
One can see this approach as a two-step dimension reduction 
procedure, where first the structure model has been a priori reduced and afterwards the dimension
reduction for the fluid part is applied.

The structure domain will be $S_h = \omega\subset\R^{2}$
(independent of $h$), and we assume that the plate is linear, visco-elastic and in the bending regime, 
i.e.~the elasticity operator $\el$ from (\ref{1.eq:EL}) is given by
\begin{equation}\label{def:L_e+0}
\el(\bu,\pa_t\bu) = B\Delta'^4\eta + \vartheta\,\Delta'^4\pa_t\eta\,,
\end{equation}
where the structure displacement is assumed to be of the form 
$\bu(\bx,t) = (\bs 0,\eta(x',t))\in\R^3$, and parameters $B$ and $\vartheta$ describe
material properties: rigidity and visco-elasticity, respectively.
Since the structure model is considered as a boundary condition on the top boundary of the fluid domain,
the dynamic interface condition becomes the balance of forces 
on the top boundary of the fluid domain which is achieved by adding the right hand side 
$\bs g = -\sigma_f(\bv,p)\bs e_3\cdot\bs e_3$ in (\ref{1.eq:EL}).
The kinematic condition reads $\pa_t\bu = \bv$.
Analogous 2D model has been investigated in \cite{PaSt06}.
A similar model with 
\begin{equation}\label{def:L_e+0adrijana}
\el(\bu,\pa_t\bu) = B\Delta'^4\eta - M\,\Delta'^2\eta\,,
\end{equation}
where $M>0$ accounts for contribution of the horizontal tension to the vertical displacement,
also in 2D, was discussed in \cite{CuMP18} and analogous results to those
in \cite{PaSt06} were obtained.

In the following we work with (\ref{def:L_e+0}). Since horizontal components of the fluid velocity
are zero on $\omega$ we have the Korn equality for $\bv$. Hence, testing equations (\ref{1.eq:St}) 
and (\ref{1.eq:EL}) with assumed smooth solutions $\bv$ and $\pa_t\eta$,
respectively, yields to the following energy estimate: for every 
$t>0$
\begin{align}
\frac{\vro_f}{2}\int_{\Omega_\eps}|\bv(t)|^2 \dd \bx
&+ \frac{\nu}{2}\int_0^t\!\!\int_{\Omega_\eps}|\nabla\bv(s)|^2 \dd\bx \dd s 
+ \frac{\vro_s}{2}\int_{\omega}|\pa_t\eta(t)|^2\dd x' \nonumber\\
&+ \frac{B}{2} \int_{\omega}|\Delta'\eta(t)|^2\dd x' 
+ \vartheta\int_0^t\!\!\int_{\omega}|\Delta'\pa_t\eta(s)|^2\dd x'\dd s \label{eq.BE:e+0}
\leq Ct\eps^3\,.
\end{align}
The obtained energy estimate motivates the following structure function space to be specified
\begin{align*}
\V_S(0,T;\omega) = W^{1,\infty}(0,T;L^2(\omega))\cap H^1(0,T;V_S(\omega))\,,
\end{align*}
where $V_S(\omega) = \left\{\eta\in H^2(\omega)\ :\ 
\eta \text{ is }\omega\text{-periodic} \right\}$, and
the solution space of the coupled problem (\ref{1.eq:St})-(\ref{def:IC}) 
with (\ref{def:L_e+0}) is then given by
\begin{align}\label{2.def:sol_space}
\V(0,T;\Omega) = \big\{(\bv,\eta)\in \V_F(0,T;\Omega_\eps)&\times \V_S(0,T;\omega)\ :\\
\bv(t) &= (\bs 0,\pa_t\eta(t))\text{ on }\omega\text{ for a.e. }t\in(0,T)\big\}.\nonumber
\end{align}

The wellposedness of the problem (\ref{1.eq:St})-(\ref{def:IC}) in two space dimensions
with the structure operator (\ref{def:L_e+0}) was addressed in \cite{PaSt06} and 
the following regularity of the weak solution $(\bv^\eps,\eta^\eps)$ is obtained:
\begin{align*}
\pa_t \bv^\eps\in L^2(0,T;L^2(\Omega_{\eps};\R^2))\quad\text{and}\quad
\pa_{tt}\eta^\eps\in L^2(\omega\times (0,T))\,,&\quad\text{(time regularity)}\,,\\
\bv^\eps\in L^\infty(0,T;H^1(\Omega_\eps;\R^2))\quad\text{and}\quad 
\eta^\eps\in L^\infty(0,T;H^2(\omega))\,,&\quad\text{(space regularity)}\,.
\end{align*}
Moreover, there exists a unique pressure $p^\eps\in L^2(0,T;H^1(\Omega_\eps))$ such that 
$(\bv^\eps,p^\eps,\eta^\eps)$ solves (\ref{1.eq:St})-(\ref{1.eq:EL}) in the classical sense.
Even tough the result in \cite{PaSt06} is stated and proved in 2D case, 
the Galerkin construction scheme and a priori estimates would provide
the same results also in the 3D case. Here we work with the 3D case.


The aim is here to obtain a nontrivial limit behavior of the original system as the small parameter $\eps$ 
tends to zero. The same problem has been analyzed in \cite{PaSt06}, but using the asymptotic
expansion techniques, like also in \cite{CuMP18}. Here we follow another concept developed in 
\cite{BuMu19} for analysis of an $\eps+h$ problem, which is based on careful quantitative energy
estimates. For that purpose the following scaling ansatz is assumed
\begin{enumerate}
  \item[(S1)] $B = \hat{B}\eps^{-\kappa}$ and $\vro_s = \hat{\vro}_s\eps^{-\kappa}$ for some
$\kappa>0$ and $\hat{B},\hat{\vro}_s>0$ independent of $\eps$;
  \item[(S2)] $\T = \eps^\tau$ for some $\tau\in\R$.
\end{enumerate}
Scaling (S1) takes into account large rigidity of the structure where
$\kappa$ may be interpreted as a measure of the structure rigidity \cite{Cia88}, and
(S2) is the choice of the time scale $\T$ depending on $\eps$. For now the point of the above
scalings is purely mathematical with aim of finding a relation between free parameters $\kappa$
and $\tau$ which will ensure the nontrivial coupled behavior of the system in the reduced model.

Next we perform the geometric change of variables from the thin fluid domain $\Omega_\eps$ to the 
reference domain $\Omega_- = (0,1)^2\times(-1,0)$ and obtain the uniform energy estimates on $\Omega_-$. 
Let us denote by $\bv(\eps)$ and $\eta(\eps)$ weak solutions to the rescaled system and denote the scaled
differential operators by $\nabla_\eps\bv = (\nabla'\bv,\eps^{-1}\pa_3\bv)$ (the scaled gradient) and
$\diver_\eps\bs v = \pa_1 v_1 + \pa_2 v_2 + \frac{1}{\eps} \pa_3v_3$ (the scaled divergence operator).
The energy estimate (\ref{eq.BE:e+0}) on the 
reference domain and in the rescaled time then reads: for a.e.~$t\in(0,T)$ it holds
\begin{align}\label{ineq:e_ref_e+0}
\frac{\vro_f}{2}\eps\int_{\Omega_-} |\bv(\eps)(t)|^2 \dd\by
 &+ \frac{\nu \eps^{\tau+1}}{2}\int_0^t\!\!\int_{\Omega_-}|\nabla_\eps\bv(\eps)|^2 \dd\by\dd s +
\frac{\vro_s}{2}\eps^{-2\tau-\kappa}\int_{\omega}|\pa_t\eta(\eps)(t)|^2\dd x' \nonumber\\
&+ \frac{B}{2}\eps^{-\kappa} \int_{\omega}|\Delta'\eta(\eps)(t)|^2\dd x' 
+ \vartheta\eps^{-\tau}\int_0^t\!\!\int_{\omega}|\Delta'\pa_t\eta(\eps)(s)|^2\dd x'\dd s \leq C \eps^{\tau+3}\,,
\end{align}
where $'$ denotes horizontal variables and respective operators.
Estimate (\ref{ineq:e_ref_e+0}) provides the uniform bound for the fluid velocity
\begin{equation*}
\int_0^T\!\!\int_{\Omega_-}|\pa_3\bv(\eps)|^2 \dd\by\dd s \leq C\eps^4\,,
\end{equation*}
which, using the no-slip boundary condition on the bottom of the channel, directly implies
\begin{equation*}
\|\bv(\eps)\|_{L^2(0,T;L^2(\Omega_-)} \leq C\eps^2\,.
\end{equation*}
This motivates to rescale the fluid velocity according to
$\bar\bv(\eps) = \eps^{-2}\bv(\eps)$ 
and neglecting the bar notation, uniform a priori estimates imply 
the following weak convergence results (on a subsequence as $\eps\downarrow0$):
\begin{align}\label{eq:vel_conv}
\bv(\eps) \rightharpoonup \bv\quad \text{and}\quad 
\pa_3\bv(\eps) \rightharpoonup \pa_3\bv\quad\text{weakly in }L^2(0,T;L^2(\Omega_-;\R^3))\,.
\end{align}

The pressure in the system is treated in a standard manner.
Define $\pi_{\eps}(t)=\int_{\Omega_{-}}p(\eps)(\by,t)\dd\by$ to be the mean value 
of the pressure at time $t\in(0,T)$. 
First, the zero mean value part of the pressure $p(\eps)-\pi_\eps$ is estimated 
in a classical way by utilizing the Bogovski operator for the construction of an appropriate 
test function (cf.~\cite[Section 3.1]{BuMu19}), which provides
\begin{equation*}
\|p(\eps) - \pi_\eps\|_{L^2(0,T;L^2(\Omega_-))} \leq C\,.
\end{equation*}
Estimating the mean value as $\|\pi_\eps\|_{L^2(0,T)}\leq C$ (cf.~again \cite[Section 3.1]{BuMu19})
yields the uniform bound for the pressure 
$\|p(\eps)\|_{L^2(0,T;L^2(\Omega_-))}\leq C$
which implies the existence of $p\in L^2(0,T;L^2(\Omega_-))$ such that
(on a subsequence as $\eps\downarrow0$)
\begin{align}\label{eq:pressure_conv}
p(\eps) \rightharpoonup p\quad \text{weakly in }L^2(0,T;L^2(\Omega_-))\,.
\end{align}

According to (\ref{ineq:e_ref_e+0}), for the structure displacement we have the bound
\begin{equation*}
\esssup_{t\in(0,T)}\int_{\omega}|\Delta'\eta(\eps)(t)|^2\dd x' \leq C\eps^{\tau + \kappa + 3}\,,
\end{equation*}
which due to the Poincar\'e inequality, since $\int_\omega\eta(\eps)(t)\dd x' = 0$ for a.e.~$t\in(0,T)$,
yields
\begin{equation}
\|\eta(\eps)\|_{L^\infty(0,T;H^2(\omega))} \leq C\eps^{(\tau + \kappa + 3)/2}\,.
\end{equation}
Rescaling $\eta(\eps)$ according to $\bar{\eta}(\eps) = \eps^{-(\tau + \kappa + 3)/2}\eta(\eps)$
and taking all previous rescalings into account yields the weak form on the reference domain
which includes the pressure:
\begin{align}
\nonumber
-\vro_f \eps^{3-\tau} \int_0^T\!\!\int_{\Omega_-}\bv(\eps)\cdot\pa_t\bs\phi\,\dd \by \dd t 
+ 2\nu\eps^3\int_0^T\!\!\int_{\Omega_-}\sym\nabla_\eps\bv(\eps):\sym\nabla_\eps\bs\phi\,\dd \by\dd t \\
- \eps\int_0^T\!\!\int_{\Omega_-}p(\eps)\diver_\eps\bs\phi\,\dd \by\dd t + \label{eq:weak_res_e+0}
\vro_s\eps^{\delta - 2\tau}\int_0^T\!\! \int_{\omega}\eta(\eps)\pa_{tt}\psi \,\dd x'\dd t \\  + 
\eps^\delta\int_0^T\!\!\int_{\omega}\left(B\Delta'\eta(\eps)\Delta'\psi \nonumber
- \vartheta\eps^{\kappa-\tau}\Delta'\eta(\eps)\Delta'\pa_t\psi\right)\dd x'\dd t \\
 = \eps\int_0^T\!\!\int_{\Omega_-}\bs f(\eps)\cdot\bs\phi \,\dd \by\dd t\,, \nonumber
\end{align}
for all $(\bs\phi,\psi)\in C_c^2\left([0,T); V(\Omega_-)\times V_S(\omega)\right)$ such that 
 $\bs\phi(t) = (\bs0,\psi(t)) \text{ on }\omega \text{ for all }t\in [0,T)$,
 and where $\delta = (\tau + \kappa + 3)/2 - \kappa$.
Omitting the divergence free condition in fluid test functions, the above fluid space is  
$V(\Omega_-) = \left\{\bv\in H^1(\Omega_-;\R^3)\ 
:\ \left.\bv\right|_{\{y_3=-1\}}=0\,,\ \bv \text{ is }\omega\text{-periodic} \right\}$.
 
In order to realize a nontrivial coupling between the fluid and the structure part 
in the reduced model we need to adjust $\delta=0$. 
Namely, in this case the fluid pressure will balance the structure
bending. This condition then yields the choice of the right time scale $\T = \eps^\tau$ with 
\begin{equation}\label{eq:rel_tau_e+0}
\tau = \kappa - 3\,.
\end{equation} 
\begin{remark}
The obtained relation (\ref{eq:rel_tau_e+0}) relates the rigidity of the structure 
with the corresponding time scale
which ``sees'' the interaction between the subsystems in the reduced model. 
(\ref{eq:rel_tau_e+0}) is consistent
with assumptions and results obtained in \cite{PaSt06} and \cite{CuMP18}. In \cite{PaSt06}, assumed
time scale corresponds to $\tau=0$ and the nontrivial coupling between the subsystems in the reduced model
is realized for the structure rigidity which corresponds to $\kappa=3$.
On the other hand, in \cite{CuMP18} the assumed
time scale is given by $\tau=-2$ and the nontrivial coupling is realized when the structure rigidity 
corresponds to $\kappa=1$.
\end{remark}
The leading order terms in (\ref{eq:weak_res_e+0}) ($O(1)$ with respect to $\eps$) are the pressure
term in the fluid part and the bending term in the structure.
Hence, under additional assumption $\tau < 0$, which ensures that the inertial term of the structure
vanishes,
the limit form of (\ref{eq:weak_res_e+0}) (on a subsequence as $\eps\downarrow0$) reads
\begin{align}\label{eq:limit_fsi_e+0}
-\int_0^T\!\!\int_{\Omega_-}p\pa_3\phi_3\,\dd \by\dd t 
+ B\int_0^T\!\!\int_{\omega} \Delta'\eta\Delta'\psi\, \dd x'\dd t & = 0\,.
\end{align}
Taking a test function $\bs\phi$ in (\ref{eq:limit_fsi_e+0}) which is compactly supported in space
and taking $\psi = 0$, it follows from (\ref{eq:limit_fsi_e+0}) that limit pressure is independent
of the vertical variable $z_3$, and therefore $p$ (although $L^2$-function)
has the trace on $\omega$. Since $\phi_3 = \psi_3$ on $\omega\times(0,T)$, after integrating by parts 
in the pressure term, the limit form (\ref{eq:limit_fsi_e+0}) then becomes
\begin{align}\label{eq:vert_e+0}
-\int_0^T\!\!\int_{\omega}p\, \psi\,\dd x'\dd t
+ B\int_0^T\!\!\int_{\omega} \Delta'\eta\Delta'\psi\, \dd x'\dd t & = 0\,
\end{align}
for arbitrary $\psi\in C_c^2([0,T);\Hper^2(\omega))$.

In order to close the limit model, we need to further explore on the fluid part.
First, multiplying $\diver_\eps\bv(\eps) = 0$ by a test function
$\varphi\in C_c^1([0,T);\Hper^1(\omega))$, integrating over space and time, integrating by parts and
employing the rescaled kinematic condition for the vertical component $v_d(\eps)/\eps = \pa_t \eta(\eps)$ 
a.e.~on $\omega\times(0,T)$,
we find
\begin{align}\label{eq:divh}
-\int_0^T\!\!\int_{\Omega_-}(v_1(\eps)\pa_1\varphi + v_2(\eps)\pa_2\varphi)\, \dd \by\dd t -
\int_0^T\!\!\int_{\omega}\eta(\eps)\pa_t\varphi\,\dd y'\dd t = 0\,,
\end{align}
which (on a subsequence as $\eps\downarrow0$) implies
\begin{align}\label{eq:div_limit}
-\int_0^T\!\!\int_{\Omega_-} (v_1\pa_1\varphi + v_2\pa_2\varphi)\,\dd \by\dd t - 
\int_0^T\!\!\int_{\omega}\eta\pa_t\varphi\,\dd y'\dd t = 0\,,
\end{align}
for all $\varphi\in C_c^1([0,T);\Hper^1(\omega))$.
This relates the limit vertical displacement of the structure with limit horizontal fluid velocities.
For the vertical fluid velocity, the divergence free equation and the no-slip boundary
condition imply $v_3=0$.
Finally, to close the reduced model, relation between horizontal fluid velocities $v_\alpha$ and 
the pressure $p$ is obtained from (\ref{eq:weak_rescaled}) by appropriate choice of test functions: 
$\bs\phi = (\phi_1/\eps,\phi_2/\eps,0)$ with 
$\phi_\alpha\in C_c^1([0,T);C_c^\infty(\Omega_-))$ and
$\psi = 0$. Convergence results (\ref{eq:vel_conv}) and (\ref{eq:pressure_conv}) then yield the 
limit equation
\begin{align}\label{eq:pv}
\nu\int_0^T\!\!\int_{\Omega_-}(\pa_3v_1\pa_3\phi_1 + \pa_3v_2\pa_3\phi_2)\,\dd\by\dd t
&- \int_0^T\!\!\int_{\Omega_-} p(\pa_1\phi_1 + \pa_2\phi_2)\,\dd\by\dd t \\
 &\quad= \int_0^T\!\!\int_{\Omega_-}(f_1\phi_1 + f_2\phi_2) \,\dd\by\dd t\,. \nonumber
\end{align}

Since the pressure $p$ is independent of the vertical variable $y_3$, equation (\ref{eq:pv})
can be solved for $v_\alpha$ explicitly in terms of $y_3$ and $p$. 
The boundary conditions are inherited from the original no-slip conditions, 
i.e.~$v_\alpha(\cdot,-1,\cdot) = v_\alpha(\cdot,0,\cdot) = 0$.
Explicit solution of $v_\alpha$ from (\ref{eq:pv}) is then given by
\begin{equation}\label{eq:v_a}
v_\alpha(\by,t) = \frac{1}{2\nu}y_3(y_3+1)\pa_\alpha p(y',t) + F_\alpha(\by,t)\,,
\quad (\by,t)\in\Omega_-\times(0,T)\,,\ \alpha=1,2\,,
\end{equation} 
where $\displaystyle F_\alpha(\cdot,y_3,\cdot) 
= \frac{y_3+1}{\nu}\int_{-1}^{0} \zeta_3 f_\alpha(\cdot,\zeta_3,\cdot)\,\dd \zeta_3 + 
\frac{1}{\nu}\int_{-1}^{y_3}(y_3-\zeta_3) f_\alpha(\cdot,\zeta_3,\cdot)\,\dd \zeta_3$. 
Replacing $v_\alpha$ from (\ref{eq:v_a}) into equation (\ref{eq:div_limit}) we obtain 
a Reynolds type equation
\begin{equation}\label{eq:reynolds2}
\int_0^T\!\!\int_{\omega}\left(-\frac{1}{12\nu}\Delta'p - F
+ \pa_t \eta\right)\varphi\,\dd y'\dd t = 0\,,
\end{equation}
where $\displaystyle F(y',t) = -\int_{-1}^0 (\pa_1 F_1 + \pa_2 F_2)\,\dd y_3$.
Combining the latter with equation (\ref{eq:vert_e+0})
we finally obtain the reduced model in terms of the sixth-order evolution equation
for the vertical displacement  
\begin{equation}\label{eq:eta_evol}
\pa_t \eta 
- \frac{B}{12\nu}(\Delta')^3\eta = F \,. 
\end{equation}
Equation (\ref{eq:eta_evol}) is 
accompanied by trivial initial data $\eta(0) = 0$ and periodic boundary conditions.

Based on the reduced model, i.e.~knowing $\eta$ solely, we are able to recover approximate
solutions to the original FSI problem. The limit pressure $p$ and horizontal fluid velocities $v_\alpha$ 
are calculated accor\-ding to (\ref{eq:vert_e+0}) and (\ref{eq:v_a}), respectively.  
The approximate fluid velocity is then defined by
\begin{align}\label{eq:Approx}
{\bs\va}^\eps(\bx,t)=\eps^2\Big(v_1(x',\frac{x_3}{\eps},t),v_2(x',\frac{x_3}{\eps},t),\va^\eps_3\Big)\,,
\quad (\bx,t)\in\Omega_\eps\times(0,T)\,,
\end{align}
where 
$\displaystyle\va_3^\eps(\bx,t)=-\eps\int_{-1}^{x_3/\eps}(\pa_1v_1
+ \pa_2v_2)(x',\xi,t)\,\dd \xi\,,$
and the approximate pressure by
$\paa^\eps(\bx,t) = p(x',t)$\,for all $(\bx,t)\in\Omega_\eps\times(0,T)$.
Moreover, the approximate vertical displacement of the structure is defined by  
\begin{equation}\label{def:aeta_eps}
\upeta^\eps(x',t) = \eps^\kappa \eta(x',t)\,,\quad (x',t)\in\omega\times(0,T)\,.
\end{equation}
The following theorem is the key result of this section.
\begin{theorem}\label{tm:EE-e+0} 
Let $(\bv^\eps,p^\eps,\eta^\eps)$ be the classical solution to the FSI problem 
(\ref{1.eq:St})-(\ref{def:IC}) with (\ref{def:L_e+0}) in rescaled time,  
let $(\bs\va^\eps,\paa^\eps,\upeta^\eps)$ be approximate solution constructed from the reduced model
as above, and assume that $0 < \kappa\leq 5/2$,
then
\begin{align*}
\|\bv^\eps - \bs \va^\eps\|_{L^2(0,T;L^2(\Omega_\eps))} &\leq C\eps^{3}\,,\\
\|p^\eps - \paa^\eps\|_{L^2(0,T;L^2(\Omega_\eps))} &\leq C\eps\,,\\
\|\eta^\eps - \upeta^\eps\|_{L^\infty(0,T;H^2(\omega))} &\leq C\eps^{\kappa+1/2}\,,
\end{align*}
where $C>0$ denotes a generic positive constant independent of $\eps$.
\end{theorem}

\begin{remark}
Error estimates for approximate solutions of the same problem have been derived 
in \cite[Theorem 6.1]{PaSt06}. However, employing the asymptotic expansion techniques
will provide good error estimates only for higher-order approximations, while 
here we present optimal estimates for the zero-order approximation.
\end{remark}

\noindent{\em Proof.}
Based on the limit equation (\ref{eq:pv}), the approximate fluid velocity ${\bs\va}^\eps$ 
and the pressure $\paa^\eps$ satisfy the modified Stokes system
\begin{equation}\label{eq:modifStokes}
\varrho_f\T^{-1}\partial_t \bs\va^\eps - \diver\sigma_f(\bs\va^\eps,\paa^\eps) 
= \bs f^\eps - f_3^\eps\bs e_3 + \bs \res^\eps_f\,,
\end{equation}
where the residual term $\bs \res^\eps_f$ is given by
$\bs \res^\eps_f = \varrho_f\T^{-1}\partial_t \bs\va^\eps 
- \nu\Delta'{\bs\va}^\eps
-\nu\partial_{33} \va_3^\eps\,\bs e_3\,$ and enjoys the uniform bound
$\|\bs \res_f^\eps\|_{L^2(0,T;L^2(\Omega_{\eps}))}\leq C\eps^{3/2}\,.$
Multiplying equation (\ref{eq:modifStokes}) by a test function $\bs\phi\in C^1_c([0,T);V_F(\Omega_\eps))$,
 and then integrating over $\Omega_\eps\times(0,T)$, we find
\begin{align}\label{eq:modifStokes_weak}
-\vro_f\int_0^T\!\! \int_{\Omega_\eps}\bs\va^\eps\cdot\pa_t\bs\phi\,\dd \bx\dd t
+& 2\nu\T\int_0^T\!\!\int_{\Omega_\eps}\sym\nabla\bs\va^\eps:\sym\nabla\bs\phi\,\dd \bx\dd t \\
-\T\int_0^T\!\!\int_{\omega}\sigma_f(\bs\va^\eps,\paa^\eps)\bs\phi\cdot \bs e_3\,\dd \bx\dd t 
&= \T\int_0^T\!\!\int_{\Omega_\eps}(f_1^\eps\phi_1 + f_2^\eps\phi_2) \,\dd \bx\dd t 
+ \T\int_0^T\!\!\int_{\Omega_\eps}\bs\res_f^\eps\cdot\bs\phi \,\dd \bx\dd t \,.\nonumber
\end{align}
Expanding the boundary term, employing the pressure relation (\ref{eq:vert_e+0}) and utilizing the
definition the approximate displacement $\upeta^\eps$ we find
\begin{align}
-\vro_f\int_0^T\!\! \int_{\Omega_\eps}\bs\va^\eps\cdot\pa_t\bs\phi\,\dd \bx\dd t\nonumber
+ 2\nu\T\int_0^T\!\!\int_{\Omega_\eps}\sym\nabla\bs\va^\eps:\sym\nabla\bs\phi\,\dd \bx\dd t \\
-\vro_s^\eps\T^{-1}\int_0^T\!\!\int_{\omega} \pa_t\upeta^\eps\pa_t\psi\,\dd x'\dd t \nonumber
+ B^\eps\T\int_0^T\!\!\int_{\omega}\Delta'\upeta^\eps\Delta'\psi\,\dd x'\dd t + 
\vartheta \int_0^T\!\!\int_{\omega}\Delta'\pa_t\upeta^\eps\Delta'\psi\,\dd x'\dd t\\
= \T\int_0^T\!\!\int_{\Omega_\eps}(f_1^\eps\phi_1 + f_2^\eps\phi_2) \,\dd \bx\dd t 
+ \T\int_0^T\!\!\int_{\Omega_\eps}\bs\res_f^\eps\cdot\bs\phi \,\dd \bx\dd t 
+ \langle\bs\res_s^\eps,\psi\rangle \,,\nonumber
\end{align}
where the structure residual term is given by
\begin{equation*}
\langle\bs\res_s^\eps,\psi\rangle = 
-\vro_s^\eps\T^{-1}\int_0^T\!\!\int_{\omega} \pa_t\upeta^\eps\pa_t\psi\,\dd x'\dd t +
\vartheta \int_0^T\!\!\int_{\omega}\Delta'\pa_t\upeta^\eps\Delta'\psi\,\dd x'\dd t\,,
\end{equation*}
while $B^\eps = B\eps^{-\kappa}$ and $\vro^\eps = \vro\eps^{-\kappa}$ 
as in the scaling ansatz (S1).
Introducing the error functions $\bs e_f^\eps := \bv^\eps - \bs\va^\eps$ 
and $e_\eta^\eps:=\eta^\eps - \upeta^\eps$ we arrive to the error equation in the weak form
\begin{align}
-\vro_f\int_0^T\!\! \int_{\Omega_\eps}\bs e_f^\eps\cdot\pa_t\bs\phi\,\dd \bx\dd t\nonumber
+ 2\nu\T\int_0^T\!\!\int_{\Omega_\eps}\sym\nabla\bs e_f^\eps:\sym\nabla\bs\phi\,\dd \bx\dd t \\
-\vro_s^\eps\T^{-1}\int_0^T\!\!\int_{\omega} \pa_t e_\eta^\eps\pa_t\psi\,\dd x'\dd t \label{eq:weak_err_e+0}
+ B^\eps\T\int_0^T\!\!\int_{\omega}\Delta'e_\eta^\eps\Delta'\psi\,\dd x'\dd t + 
\vartheta \int_0^T\!\!\int_{\omega}\Delta'\pa_te_\eta^\eps\Delta'\psi\,\dd x'\dd t\\
= \T\int_0^T\!\!\int_{\Omega_\eps}f_3^\eps\phi_3 \,\dd \bx\dd t 
- \T\int_0^T\!\!\int_{\Omega_\eps}\bs\res_f^\eps\cdot\bs\phi \,\dd \bx\dd t 
- \langle\bs\res_s^\eps,\psi\rangle \,.\nonumber
\end{align}
Observe that by construction, approximate solutions satisfy the same boundary conditions as the
true solutions to the original problem. Namely, $\left.\va_\alpha^\eps\right|_{\{y_3=-1\}} = \left.\va_\alpha^\eps\right|_\omega=0$,
while for the vertical component we have
\begin{equation*}
\left.\va_3^\eps\right|_\omega = -\eps^3\int_{-1}^0(\pa_1  v_1 + \pa_2  v_2)\dd y_3 
= \eps^3\pa_t\eta = \T^{-1}\pa_t\upeta^\eps\,.
\end{equation*}
Thus, utilizing the Korn equality and 
$(\bs\phi,\psi) = (\bs e_f^\eps, \T^{-1}\pa_t e_\eta^\eps)$ as test functions in
(\ref{eq:weak_err_e+0}) we find: for a.e.~$t\in(0,T)$
\begin{align*}
\frac{\vro_f}{2}\int_{\Omega_\eps}|\bs e_f^\eps(t)|^2\dd\bx &+ \nonumber
\nu\T\int_0^t\!\!\int_{\Omega_\eps}|\nabla\bs e_f^\eps|^2\,\dd \bx\dd s
+\frac{\vro_s^\eps\T^{-2}}{2}\int_{\omega}|\pa_t e_{\eta}^\eps(t)|^2 \dd x' \\ 
& + \frac{B^\eps}{2}\int_{\omega}|\Delta' e_\eta^\eps(t)|^2\dd x' + 
\vartheta\T^{-1}\int_0^t\!\!\int_{\omega}|\Delta'\pa_t e_\eta^\eps(s)|^2\dd x'\dd s \\
&= \T\int_0^t\!\!\int_{\Omega_\eps}f_3^\eps e_{f,3}^\eps \,\dd \bx\dd s
- \T\int_0^t\!\!\int_{\Omega_\eps}\bs\res_f^\eps\cdot\bs e_f^\eps \,\dd \bx\dd s 
- \T^{-1}\langle\bs\res_s^\eps,\pa_t e_{\eta}^\eps\rangle\,.
\end{align*}

Let us now estimate the right hand side. Employing a higher-order energy estimate 
(cf.~\cite[Section 2.4]{BuMu19}) one can conclude 
$\|\pa_\alpha \bv^\eps\|_{L^2(0,T;L^2(\Omega_\eps))}\leq C\eps^{5/2}$, which combined with the 
divergence free condition on the thin domain provides 
$\|v_3^\eps\|_{L^2(0,T;L^2(\Omega_\eps))}\leq C\eps^{7/2}$. As a consequence, we easily conclude 
$\|e_{f,3}^\eps\|_{L^2(0,T;L^2(\Omega_\eps))}\leq C\eps^{7/2}$, which in further gives
\begin{equation*}
\T\left|\int_0^t\!\!\int_{\Omega_\eps}f_3^\eps e_{f,3}^\eps \,\dd \bx\dd s\right|\leq C\T\eps^4\,.
\end{equation*}
For the fluid residual we employ Cauchy-Schwarz, Poincar\'e and the Young inequality, respectively,
and obtain
\begin{align*}
\T\left|\int_0^t\!\!\int_{\Omega_\eps}\bs\res_f^\eps\cdot\bs e_f^\eps \,\dd \bx\dd s\right|
\leq C\T\int_0^t \eps \|\bs\res_f^\eps\|_{L^2(\Omega_\eps)}\|\nabla \bs e_f^\eps\|_{L^2(\Omega_\eps)}\dd s
\leq C\T\eps^5 + \frac{\nu\T}{2}\int_0^t\!\!\int_{\Omega_\eps}|\nabla\bs e_f^\eps|^2\,\dd \bx\dd s\,.
\end{align*}
Finally, we estimate the structure residual,
\begin{align*}
\T^{-1}\left|\langle\bs\res_s^\eps,\pa_t e_{\eta}^\eps\rangle \right|
\leq \vro_s^\eps\T^{-2}\left|\int_0^t\!\!\int_{\omega} \pa_t\upeta^\eps\pa_{tt}e_\eta^\eps\,\dd x'\dd s\right|
 + \vartheta\T^{-1} \left|\int_0^t\!\!\int_{\omega}\Delta'\pa_t\upeta^\eps\Delta'\pa_te_\eta^\eps\,\dd x'\dd s\right|\,.
\end{align*}
Integrating by parts in time, the first term can be estimated as
\begin{align*}
\vro_s^\eps\T^{-2}\left|\int_0^t\!\!\int_{\omega} \pa_t\upeta^\eps\pa_{tt}e_\eta^\eps\,\dd x'\dd s\right|
\leq C\T\eps^{3-2\tau} + \frac{\vro_s^\eps\T^{-2}}{4}\int_{\omega}|\pa_t e_{\eta}^\eps(t)|^2 \dd x'
 + \frac{\vro_s^\eps\T^{-2}}{4}\int_0^t\!\!\int_{\omega}|\pa_t e_{\eta}^\eps(s)|^2 \dd x'\dd s\,,
\end{align*}
while for the second term we have 
\begin{align*}
\vartheta\T^{-1} \left|\int_0^t\!\!\int_{\omega}\Delta'\pa_t\upeta^\eps\Delta'\pa_te_\eta^\eps\,\dd x'\dd s\right|
\leq C\T\eps^6 + \frac{\vartheta\T^{-1}}{2}
\int_0^t\!\!\int_{\omega}|\Delta'\pa_t e_\eta^\eps(s)|^2\dd x'\dd s\,.
\end{align*}
Summing all up and employing the Gronwall inequality we have: for a.e.~$t\in(0,T)$
\begin{align}
\frac{\vro_f}{2}\int_{\Omega_\eps}|\bs e_f^\eps(t)|^2\dd\bx &+ \nonumber
\frac{\nu\T}{2}\int_0^t\!\!\int_{\Omega_\eps}|\nabla\bs e_f^\eps|^2\,\dd \bx\dd s
+\frac{\vro_s^\eps\T^{-2}}{4}\int_{\omega}|\pa_t e_{\eta}^\eps(t)|^2 \dd x' \\ 
& + \frac{B^\eps}{2}\int_{\omega}|\Delta' e_\eta^\eps(t)|^2\dd x' + 
\frac{\vartheta\T^{-1}}{2} \int_0^t\!\!\int_{\omega}|\Delta'\pa_t e_\eta^\eps(s)|^2\dd x'\dd s  
\label{ineq:error_est_e+0}
  \leq C\T\eps^4\,
\end{align}
provided $\tau \leq -1/2$, i.e.~$\kappa\leq 5/2$.

Combining estimate (\ref{ineq:error_est_e+0}) and the Poincar\'e inequality we find
\begin{equation}
\|\bs e_f^\eps\|_{L^2(0,T;L^2(\Omega_\eps))} \leq C\eps^3\,,
\end{equation}
which is the desired error estimate for the fluid velocity. Due to the zero mean value on $\omega$,
according to (\ref{ineq:error_est_e+0}), the displacement error can be controlled as
\begin{equation}
\|e_\eta^\eps\|_{L^\infty(0,T;H^2(\omega))}^2 \leq C\eps^{2\kappa+1}\,,
\end{equation}
which is the required estimate for structure displacement. The error estimate for the pressure follows
somewhat different approach, which we omit here and refer to \cite[Section 4.5]{BuMu19} for details. 

$\hfill\Box$

\begin{remark}
Observe that the error estimate of horizontal fluid velocities relative to the norm of velocities
is $O(\sqrt{\eps})$. The same holds true for the relative error 
estimate of the pressure and displacement errors.
\end{remark}

\subsection{Linear $\eps + h$ problem}\label{sec:e+h}
In practice we often have that both layers are relatively thin of equivalent size. 
Moreover, we expect then that the
interplay between the two thicknesses: $\eps$ of the fluid and $h$ od the structure plays the 
role in the derivation of reduced models. Unlike in the previous section, here we take a full
structure model and aim to perform a simultaneous dimension reduction.

In the FSI problem (\ref{1.eq:St})-(\ref{def:IC}), the channel is now covered by an elastic 
plate of relative height $h>0$ with material configuration $S_h = (0,1)^2\times(0,h)$
and the plate is described by the linear 3D elasticity
equations, i.e.~the elasticity operator $\el$ in (\ref{1.eq:EL}) is given by
\begin{equation}\label{def:L_e+h}
\el(\bu,\pa_t\bu) = -\diver\sigma_s(\bu)\,,
\end{equation}
where $\sigma_s(\bu)=2\mu\sym\nabla\bu  + \lambda(\diver\bu )I_3$ denotes the Cauchy stress tensor, 
$\mu$ and $\lambda$ are Lam\'e constants and $I_3$ is $3\times3$ identity matrix.
The two subsystems are coupled through the
interface conditions on the fixed interface $\omega$: 
\begin{align}\label{def:kinematic_bc}
\pa_t\bu &= \bv\,, \quad \omega\times(0, \infty)\,,
\qquad \text{(kinematic -- continuity of velocities)\,,}\\
(\sigma_f(\bv,p) - \sigma_s(\bu))\bs e_3 &=0\,,\quad \omega\times(0, \infty)\,,
\qquad \text{ (dynamic -- balance of forces)\,.} \label{def:dynamic_bc}
\end{align}
The action of the fluid on the structure and vice versa is considered through the 
dynamic coupling condition (\ref{def:dynamic_bc}) and we assume the absence of the structure volume 
forces, hence, $\bs g = 0$ in (\ref{1.eq:EL}).
Following \cite{BuMu19} the basic energy inequality (\ref{eq.BEI}) in this case reads: for every $t>0$
\begin{align}
\frac{\vro_f}{2}\int_{\Omega_\eps}|\bv(t)|^2 \dd \bx
&+ \nu\int_0^t\!\!\int_{\Omega_\eps}|\sym\nabla\bv(s)|^2 \dd\bx \dd s 
+ \frac{\vro_s}{2}\int_{S_h}|\pa_t\bu(t)|^2\dd\bx \nonumber\\
&+ \mu \int_{S_h}|\sym\nabla\bu(t)|^2\dd\bx 
+ \frac{\lambda}{2}\int_{S_h}|\diver \bu(t)|^2\dd\bx \label{eq.energy_basic}
\leq Ct\eps\,,
\end{align}
and the structure function space (\ref{def:str_space}) is now specified to be
\begin{align*}
\V_S(0,T;S_h) = W^{1,\infty}(0,T;L^2(S_h;\R^3))\cap L^\infty(0,T;V_S(S_h))\,,
\end{align*}
where $V_S(S_h) = \left\{\bu\in H^1(S_h;\R^3)\ :\ 
\bu \text{ is }\omega\text{-periodic} \right\}$. 
The solution space of the coupled problem (\ref{1.eq:St})-(\ref{def:IC}) 
with (\ref{def:L_e+h})-(\ref{def:dynamic_bc}) is now given by
\begin{align}\label{2.def:sol_space}
\V(0,T;\Omega) = \big\{(\bv,\bu)\in \V_F(0,T;\Omega_\eps)&\times \V_S(0,T;S_h)\ :\\
\bv(t) &= \pa_t\bu(t)\text{ on }\omega\text{ for a.e. }t\in(0,T)\big\}.\nonumber
\end{align}

The wellposedness of the problem (\ref{1.eq:St})-(\ref{def:IC}) 
with (\ref{def:L_e+h})-(\ref{def:dynamic_bc}) in the sense of Definition \ref{def:weak_sol} 
is addressed in \cite{BuMu19} and the regularity of the weak solution $(\bv^\eps,\bu^h)$ required
for the subsequent analysis is obtained as follows:
\begin{align*}
\pa_t \bv^\eps\in L^\infty(0,T;L^2(\Omega_{\eps};\R^3))\quad\text{and}\quad
\pa_{tt}\bu^h\in L^\infty(0,T;L^2(S_h;\R^3))\,,&\quad\text{(time regularity)}\,,\\
\bv^\eps\in L^\infty(0,T;H^2(\Omega_\eps;\R^3))\quad\text{and}\quad 
\bu^h\in L^\infty(0,T;H^2(S_h;\R^3))\,,&\quad\text{(space regularity)}\,.
\end{align*}
Moreover, a refined energy estimate is obtained (cf.~\cite[Section 2.4]{BuMu19}): 
for a.e.~$t\in[0,T)$
\begin{align}\label{ineq:energEH}
\frac{\vro_f}{2}\int_{\Omega_\eps}|\bv^\eps(t)|^2 \dd\bx
 & + \frac{\nu}{2}\int_0^t\!\!\int_{\Omega_\eps}|\nabla\bv^\eps|^2 \dd\bx\dd s
 + \frac{\vro_s}{2}\int_{S_h}|\partial_t\bu^h(t)|^2 \dd\bx \\
& + \int_{S_h}\Big(\mu|\sym\nabla\bu^h(t)|^2 \nonumber
+ \frac{\lambda}{2}|\diver\bu^h(t)|^2\Big)\dd\bx \leq C t \eps^{3}\,,
\end{align}
which will be the cornerstone for the derivation of the reduced model.
Furthermore, there exists a unique pressure $p^\eps\in L^2(0,T;H^1(\Omega_\eps))$ such that 
$(\bv^\eps,p^\eps,\bu^h)$ solves (\ref{1.eq:St})-(\ref{1.eq:EL}) in the classical sense.

A similar problem has been analyzed in \cite{PaSt14}, where starting from a 2D/2D linear FSI problem of
type $1+h$, the FSI problem analyzed in \cite{PaSt06} has been justified.
Our main aim here is to obtain
a nontrivial limit behavior of the original system as both small parameters $\eps$ and $h$ 
simultaneously vanish. In order to achieve that,
we need to assume some scaling ansatz:
\begin{enumerate}
  \item[(S1)] $\eps = h^\gamma$ for some $\gamma>0$ independent of $h$;
  \item[(S2)] $\mu = \hat{\mu}h^{-\kappa}$, $\lambda = \hat{\lambda}h^{-\kappa}$ and
$\vro_s = \hat{\vro}_sh^{-\kappa}$ for some
$\kappa>0$ and $\hat{\mu},\ \hat{\lambda}$, $\hat\vro_s$ independent of $h$;
  \item[(S3)] $\T = h^\tau$ for some $\tau\in\R$.
\end{enumerate}
Scaling (S1) is a geometric relation between small parameters, (S2) like in the
previous section, takes into account the large rigidity of the structure where
$\kappa$ may be interpreted as a measure of the structure rigidity \cite{Cia88}, and
(S3) simply means the choice of the time scale $\T$ depending on $h$. Again, the point of the above
scalings is for now purely mathematical and we aim to find a relation between free parameters 
$\gamma$, $\kappa$ and $\tau$, similar to (\ref{eq:rel_tau_e+0}), which will ensure the 
nontrivial limit behavior of the full system. 

Performing the standard change of variables we move to the reference domain 
$\Omega_-\cup\omega\cup S_+$ and obtain the uniform energy estimates there. 
Let us denote by $\bv(\eps)$ and $\bu(h)$ weak solutions to the rescaled system and
by $\nabla_\eps$ and $\nabla_h$ the corresponding scaled gradients, then 
the energy estimate (\ref{ineq:energEH}) on the 
reference domain and in rescaled time reads: for a.e.~$t\in(0,T)$ it holds
\begin{align}\label{ineq:energy_ref}
\frac{\vro_f}{2}\eps\int_{\Omega_-}& |\bv(\eps)(t)|^2 \dd\by
 + \frac{\nu h^\tau}{2}\eps\int_0^t\!\!\int_{\Omega_-}|\nabla_\eps\bv(\eps)|^2 \dd\by\dd s 
 + \frac{\vro_s}{2}h^{-\kappa-2\tau+1}\int_{S_+}|\partial_t\bu(h)(t)|^2 \dd\bz \\
& + h^{-\kappa+1}\int_{S_+}\Big(\mu|\sym\nabla_h\bu(h)(t)|^2 \nonumber
+ \frac{\lambda}{2}|\diver_h\bu(h)(t)|^2\Big)\dd\bz \leq Ch^\tau \eps^{3}\,.
\end{align}

The rescaled energy estimate (\ref{ineq:energy_ref}) gives us for the fluid part same uniform bounds 
and the same convergence results as in (\ref{eq:vel_conv}) and (\ref{eq:pressure_conv}).
For the structure part, the energy estimate (\ref{ineq:energy_ref}) provides an $L^\infty$-$L^2$ estimate of 
the symmetrized scaled gradient of the displacement,
\begin{equation}
\esssup_{t\in(0,T)}\int_{S_+}|\sym\nabla_h\bu(h)|^2\dd\bz\leq {Ch^{3\gamma-1+\tau + \kappa}}\,,
\end{equation}
which motivates rescaling of the structure displacements according to
$\bar\bu(h) = {h^{-(3\gamma-1+\tau + \kappa)/2}}\bu(h)$. 

The uniform bound on the symmetrized scaled gradient of displacements motivates to invoke
the framework of the Griso decomposition \cite{Gri05} --- 
for every $h>0$, scaled structure displacement $\bu(h)$ is, at almost every time instance $t\in(0,T)$, 
decomposed into a sum of so called elementary plate 
displacement and warping: 
\begin{align}\label{eq:grisodec}
\bu(h)(\bz) = \bs w(h)(z') + 
\bs r(h)(z')\times (z_3-\frac12)\bs e_3 + \tilde{\bu}(h)\,,
\end{align}
where
\begin{align*}
\bs w(h)(t,z') = \int_{0}^1\bu(h)(t,\bz)\dd z_3\,,\quad 
\bs r(h)(t,z') = \frac{3}{h}\int_0^1(z_3-\frac12)\bs e_3\times \bu(h)(t,\bz)\dd z_3\,,
\end{align*}
$\tilde{\bu}(h)\in L^\infty(0,T;H^1(S_+))$ is the warping term, and $\times$ denotes 
the cross product in $\R^3$. Moreover, the following uniform estimate holds
\begin{align*}
\|\sym\nabla_h \left(\bs w(h)(z') 
+ \bs r(h)(z')\times (z_3-1/2)\bs e_3\right)\|^2_{L^\infty(0,T;L^2(S_+))} + 
\|\nabla_h\tilde{\bu}(h)\|^2_{L^\infty(0,T;L^2(S_+))} \\
+ \frac{1}{h^2}\|\tilde{\bu}(h)\|^2_{L^\infty(0,T;L^2(S_+))} & \leq C\,,
\end{align*}
with $C>0$ independent of $h$ and $\bu(h)$.
Following \cite[Theorem 2.6]{Gri05}, the above
uniform estimate implies the existence of 
a sequence of in-plane translations $\bs a(h) = (a_1(h),a_2(h))\subset (L^\infty(0,T))^2$, as well
as limit displacements $\eta_1, \eta_2\in L^\infty(0,T;\Hper^1(\omega))$, $\eta_3\in L^\infty(0,T;\Hper^2(\omega))$
and 
$\bar \bu\in L^2(\omega;H^1((0,1);\R^3))$
such that
the following weak-$\star$ convergence results hold:
\begin{align}
w_\alpha(h)- a_\alpha(h) &\overset{\ast}{\rightharpoonup} \eta_\alpha \quad\text{in }
\ L^\infty(0,T;\Hper^1(\omega))\,,\quad \alpha = 1,2\,,\label{eq:w1_conv} \\
hw_3(h) &\overset{\ast}{\rightharpoonup} \eta_3 \quad\text{in }\ L^\infty(0,T;\Hper^1(\omega))\,,\label{eq:w3_conv} \\
u_\alpha(h) - a_\alpha(h) &\overset{\ast}{\rightharpoonup} \eta_\alpha - (z_3-\frac12)\pa_\alpha \eta_3 
\quad\text{in }\ L^\infty(0,T;H^1(S_+))\,,\quad \alpha = 1,2\,,\label{eq:lsd1}\\
hu_3(h) &\overset{\ast}{\rightharpoonup} \eta_3  \quad\text{in }\ L^\infty(0,T;H^1(S_+))\,,\label{eq:lsd3}\\
\sym\nabla_h\bu(h)  &\overset{\ast}{\rightharpoonup}  \imath\Big(\sym\nabla' (\eta_1,\eta_2)  - (z_3-\frac12)\nabla'^2\eta_3 \Big)
+ \sym \left(\bs e_3\otimes(\pa_3\bar\bu)\right). \label{eq:lss}
\end{align}
%
 

Taking all the above rescalings into account,
the weak form, which now includes the pressure, on the reference domain reads
\begin{align}
\nonumber
-\vro_f h^{-\tau} \eps^3 \int_0^T\!\!\int_{\Omega_-}\bv(\eps)\cdot\pa_t\bs\phi\,\dd \by \dd t 
+ 2\nu\eps^3\int_0^T\!\!\int_{\Omega_-}\sym\nabla_\eps\bv(\eps):\sym\nabla_\eps\bs\phi\,\dd \by\dd t \\
- \eps\int_0^T\!\!\int_{\Omega_-}p(\eps)\diver_\eps\bs\phi\,\dd \by\dd t + \label{eq:weak_rescaled}
\vro_sh^{\delta - 2\tau}\int_0^T\!\! \int_{S_+}\bu(h)\cdot\pa_{tt}\bs\psi \,\dd \bz\dd t \\  + 
h^\delta\int_0^T\!\!\int_{S_+}(2\mu\sym\nabla_h\bu(h):\sym\nabla_h\bs\psi \nonumber
+ \lambda\diver_h\bu(h)\diver_h\bs\psi)\,\dd \bz\dd t \\
 = \eps\int_0^T\!\!\int_{\Omega_-}\bs f(\eps)\cdot\bs\phi \,\dd \by\dd t\,, \nonumber
\end{align}
for all $(\bs\phi,\bs\psi)\in C_c^2\left([0,T); V(\Omega_-)\times V_S(S_+)\right)$ such that 
 $\bs\phi(t) = \bs\psi(t) \text{ on }\omega \text{ for all }t\in [0,T)$,
 and where $\delta = (3\gamma-1+\tau + \kappa)/2+1 - \kappa$.

In order to realize a nontrivial coupling in the reduced model on the limit as $h\downarrow0$,
like in the previous section we need to adjust the parameter $\delta$. 
The linear theory of plates (cf.~\cite[Section 1.10]{Cia97}) suggests $\delta=-1$. Namely, the fluid 
pressure which is here $O(1)$ is acting as a normal force on the structure and therefore 
has to balance the structure stress terms in the right way.
This condition then yields the choice of the right time scale $\T = h^\tau$ with 
\begin{equation}\label{eq:rel_tau}
\tau = \kappa - 3\gamma - 3\,.
\end{equation} 

\begin{remark}
Relation (\ref{eq:rel_tau}) is analogous to (\ref{eq:rel_tau_e+0}) and mathematically has two free
parameters: measure of the structure rigidity $\kappa$ and geometric relation between thicknesses $\gamma$,
which fix the choice of the right time scale that ``sees'' the interaction between the two
subsystems in the reduced model. 
\end{remark}

Exploring the structure of elasticity equations and taking appropriate test functions which 
imitate the shape of the limit of scaled displacements 
(\ref{eq:lsd1})--(\ref{eq:lsd3}), i.e.~$\bs\psi = (h\psi_1,h\psi_2,\psi_3)$ satisfying 
$\pa_1\psi_3 + \pa_3\psi_1 = \pa_2\psi_3 + \pa_3\psi_2 = \pa_3\psi_3 = 0$ 
(cf.~\cite[Theorem 1.4-1]{Cia97}) and
$\bs \phi = (h\phi_1,h\phi_2,\phi_3)$, under additional assumption $\tau < -1$, 
the weak limit form of (\ref{eq:weak_rescaled}) (on a subsequence as $h\downarrow0$, 
as well as $\eps\downarrow0$ ) reads
\begin{align}
-\int_0^T\!\!\int_{\Omega_-}p\pa_3\phi_3\,\dd \by\dd t 
 \label{eq:limit_fsi}
+ \int_0^T\!\!\int_{S_+}\Big(2\mu\big(\sym\nabla' (\eta_1,\eta_2) - (z_3-\frac12)\nabla'^2\eta_3\big):\sym\nabla'(\psi_1,\psi_2)\\
+ \frac{2\mu\lambda}{2\mu + \lambda}\diver\big((\eta_1,\eta_2) - (z_3-\frac12)\nabla' \eta_3 \big)\diver(\psi_1,\psi_2)\Big)\,\dd \bz\dd t & = 0\,.\nonumber
\end{align}

The obtained limit equation can be interpreted as a linear plate model \cite{Cia97} 
coupled with the limit fluid pressure acting as a normal force on the structure interface $\omega$.
Exploring in further the structure of the space of test functions for the elesticity part
\cite[Theorem 1.4-1 (c)]{Cia97}, leads to an equivalent decoupled system for 
horizontal and vertical displacements: 
%
\begin{align}\label{eq:w3_1}
-\int_0^T\!\!\int_{\omega}p\zeta_3\,\dd z'\dd t +
\int_0^T\!\!\int_{\omega}\left(\frac{4\mu}{3}\nabla'^2 \eta_3:\nabla'^2\zeta_3
+ \frac{4\mu\lambda}{3(2\mu + \lambda)}\Delta' \eta_3\Delta' \zeta_3\right) \,\dd z'\dd t 
&= 0\,
\end{align}
for arbitrary $\zeta_3\in C_c^2([0,T);\Hper^2(\omega))$ and 
\begin{align}\label{eq:w12_1}
\int_0^T\!\!\int_{\omega}\left(4\mu\sym\nabla' (\eta_1,\eta_2) :\nabla'(\zeta_1,\zeta_2)
+ \frac{4\mu\lambda}{2\mu + \lambda}\diver'\left(\eta_1,\eta_2\right)\diver'(\zeta_1,\zeta_2)\right)\dd z'\dd t
= 0\,
\end{align}
for arbitrary $\zeta_\alpha\in C_c^1([0,T);\Hper^1(\omega))$.
Equation (\ref{eq:w12_1}) implies that horizontal displacements $(\eta_1,\eta_2)$ are spatially constant
functions, and as such they will not affect the reduced model. 
Moreover, they are dominated by potentially large horizontal translations as discussed in
\cite[Section 2.5]{BuMu19}, hence we omit them in further analysis.
Thus, the limit system
(\ref{eq:limit_fsi}) is now essen\-ti\-a\-lly described with (\ref{eq:w3_1}), which relates the limit 
fluid pressure $p$ with the limit vertical displacement of the structure $\eta_3$. 

Analysis of the fluid part completely follows the lines of the previous section and results in
limit equations (\ref{eq:div_limit}) and (\ref{eq:pv}).
The limit horizontal velocities $v_\alpha$ can again be explicitly calculated in terms of $y_3$ and $p$,
but the top boundary condition for $v_\alpha$ is no longer trivial. It follows from the
interface kinematic condition that $v_\alpha(\cdot,0,\cdot) = \pa_ta_\alpha$, $\alpha = 1,2$,
where $\pa_ta_\alpha$ are translational limit velocities of the structure defined below 
and discussed in \cite[Section 3.3]{BuMu19}.
Explicit solution of $v_\alpha$ from (\ref{eq:pv}) is then given by
\begin{equation}\label{eq:v_alpha}
v_\alpha(\by,t) = \frac{1}{2\nu}y_3(y_3+1)\pa_\alpha p(y',t) + F_\alpha(\by,t) + (1+y_3)\pa_ta_\alpha\,,
\quad (\by,t)\in\Omega_-\times(0,T)\,,
\end{equation} 
where $\displaystyle F_\alpha(\cdot,y_3,\cdot) 
= \frac{y_3+1}{\nu}\int_{-1}^{0} \zeta_3 f_\alpha(\cdot,\zeta_3,\cdot)\,\dd \zeta_3 + 
\frac{1}{\nu}\int_{-1}^{y_3}(y_3-\zeta_3) f_\alpha(\cdot,\zeta_3,\cdot)\,\dd \zeta_3$. 
Replacing $v_\alpha$ from (\ref{eq:v_alpha}) into equation (\ref{eq:div_limit}) we obtain 
the same Reynolds equation (\ref{eq:reynolds2}), which
combined with equation (\ref{eq:w3_1}) results in
the reduced model in terms of the vertical displacement only.
Substituting $\eta \equiv \eta_3$ we arrive to the sixth-order evolution equation
\begin{equation}\label{eq:w3_evol}
\pa_t \eta 
- \frac{2\mu(\mu + \lambda)}{9\nu(2\mu + \lambda)}(\Delta')^3\eta = F \,, 
\end{equation}
which is (up to the coefficient in front of the spatial operator) the same as (\ref{eq:eta_evol}).

Knowing $\eta$ solely, the pressure $p$ and horizontal fluid velocities $v_\alpha$ 
are then calculated accor\-ding to
(\ref{eq:reynolds2}) and (\ref{eq:v_alpha}), respectively. Based on that we can construct
approximate fluid velocity 
\begin{align}\label{eq:Approx}
{\bs\va}^\eps(\bx,t)=\eps^2\Big(v_1(x',\frac{x_3}{\eps},t),v_2(x',\frac{x_3}{\eps},t),\va^\eps_3\Big)\,,
\quad (\bx,t)\in\Omega_\eps\times(0,T)\,,
\end{align}
where 
$\displaystyle\va_3^\eps(\bx,t)=-\eps\int_{-1}^{x_3/\eps}(\pa_1v_1
+ \pa_2v_2)(x',\xi,t)\,\dd \xi\,,$
and the approximate pressure by
$\paa^\eps(\bx,t) = p(x',t)$\,for all $(\bx,t)\in\Omega_\eps\times(0,T)$.
Moreover, the approximate displacement is defined by  
\begin{equation}\label{def:ua_h}
\bs\ua^h(\bx,t) = h^{\kappa-3}\left(h^{-\gamma}a_1 - \Big(x_3 
- \frac{h}{2}\Big)\pa_1\eta(x',t), h^{-\gamma}a_2 - \Big(x_3 - \frac{h}{2}\Big)\pa_2\eta(x',t), 
\eta(x',t) \right)\,,   
\end{equation}
for all $(\bx,t)\in\Omega_h\times(0,T)$, 
where $a_\alpha$ are horizontal time-dependent 
translations calculated by $\displaystyle a_\alpha(t) = -\int_0^t\!\!\int_\omega \pa_3F_\alpha(y',0,t)\dd y' \dd s$, 
$\alpha = 1,2$. The virtue of the reduced model is then revealed by the following convergence results
of approximate solutions.
\begin{theorem}[\cite{BuMu19}]\label{tm:EE} 
Let $(\bv^\eps,p^\eps,\bu^h)$ be the classical solution to the FSI problem 
(\ref{1.eq:St})-(\ref{def:IC}) under (\ref{def:L_e+h})-(\ref{def:dynamic_bc}) in rescaled time and  
let $(\bs\va^\eps,\paa^\eps,\bs\ua^h)$ be approximate solution constructed from the reduced model
as above. Let us additionally assume that 
$\max\{2\gamma+1, \frac74\gamma+\frac32\} \leq \kappa < 2 + 2\gamma$, then
\begin{align*}
\|\bv^\eps - \bs \va^\eps\|_{L^2(0,T;L^2(\Omega_\eps))} &\leq 
C\eps^{5/2}h^{\min\{\gamma/2,\,2\gamma-\kappa+2\}}\,,\\
\|p^\eps - \paa^\eps\|_{L^2(0,T;L^2(\Omega_\eps))} 
&\leq C\eps^{1/2}h^{\min\{\gamma/2,2\gamma-\kappa+2\}}\,,\\
\|u_\alpha^h - \ua_\alpha^h\|_{L^\infty(0,T;L^2(\Omega_h))} &\leq 
Ch^{\kappa-3/2}h^{\min\{1,\gamma/2, 2\gamma+2-\kappa\}} + 
C\sqrt{h}\|a_\alpha^h - h^{\kappa-3-\gamma}a_\alpha\|_{L^\infty(0,T)}\,,\\
\|u_3^h - \ua_3^h\|_{L^\infty(0,T;L^2(\Omega_h))} &\leq 
Ch^{\kappa-5/2}h^{\min\{1/2,\gamma/2, 2\gamma+2-\kappa\}}\,,
\end{align*}
where $C>0$ denote generic positive constants independent of $\eps$ and $h$.
\end{theorem}
\begin{remark}
Observe that, like in the previous section, the error estimate of horizontal fluid velocities 
relative to the norm of velocities
is $O(\sqrt{\eps})$ for $\kappa\leq\frac32\gamma+2$. The same holdes true for the relative error 
estimate of the pressure. For the vertical fluid velocity, which is of lower order, 
there is a lack of the error estimate. In the leading order of the structure displacement, 
namely in the vertical component, for $\kappa\leq \frac32\gamma+2$ we have the relative convergence 
rate $O(h^{\min\{1/2,\gamma/2\}})$, which means $O(\sqrt{\eps})$ for $\gamma\leq 1$ and
$O(\sqrt{h})$ for $\gamma>1$. In horizontal structure displacements, dominant part
of the error estimates are errors in horizontal translations, which are actually artifact of
periodic boundary conditions (cf.~\cite[Section 2.5]{BuMu19}). 
Neglecting these errors, which cannot be controlled in a better way, 
the relative error estimate of horizontal displacements for $\kappa\leq\frac32\gamma+2$ is 
$O(h^{\min\{1,\gamma/2\}})$, which means $O(\sqrt{\eps})$ for $\gamma\leq2$ and
$O(h)$ for $\gamma > 2$. 
Let us point out that one cannot expect better convergence rates for such first-order
approximation without dealing with boundary layers, which arrise around the interface $\omega$
due to mismatch of the interface conditions for approximate solutions.
Moreover, in \cite{MaPa01} the obtained convergence rate for the Poiseuille flow in the case 
of rigid walls of the fluid channel is $O(\sqrt{\eps})$. On the other hand, convergence rate for
the clamped Kirchhoff-Love plate is found to be $O(\sqrt{h})$ \cite{Des81}.
\end{remark}

\medskip
\noindent{\em Outline of the proof of Theorem \ref{tm:EE}.}
The proof follows the idea of the proof of Theorem \ref{tm:EE-e+0}, but mostly due to mismatch
of the interface conditions for approximate solutions, the analysis is much more involved.

Starting from equation (\ref{eq:modifStokes}) satisfied by the 
approximate fluid velocity ${\bs\va}^\eps$,
expanding the boundary terms, employing the pressure relation (\ref{eq:limit_fsi}) and utilizing the
definition the approximate displacement $\bs\ua^h$ we find
the weak form for approximate solutions to be 
\begin{align}
-\vro_f\int_0^T\!\! \int_{\Omega_\eps}\bs\va^\eps\cdot\pa_t\bs\phi\,\dd \bx\dd t \nonumber
+ 2\nu\T\int_0^T\!\!\int_{\Omega_\eps}\sym\nabla\bs\va^\eps:\sym\nabla\bs\phi\,\dd \bx\dd t\\
-\vro_s^h\T^{-1}\int_0^T\!\!\int_{S_h}\pa_t\bs\ua^h\cdot\pa_t\bs\psi\,\dd\bx\dd t 
+ \T\int_0^T\!\!\int_{S_h}\left(2\mu^h\sym\nabla\bs\ua^h :\sym\nabla\bs\psi
+ \lambda^h\diver\bs\ua^h \diver\bs\psi\right)\,\dd \bx\dd t \label{eq:weak_fsi_app}  \\
= \T\int_0^T\!\!\int_{\Omega_\eps}(f_1^\eps\phi_1 + f_2^\eps\phi_2) \,\dd \bx\dd t 
+ \T\int_0^T\!\!\int_{\Omega_\eps}\bs\res_f^\eps\cdot\bs\phi \,\dd \bx\dd t \nonumber
+ \T\int_0^T\!\!\int_{\omega}\bs\res_b^\eps\cdot\bs\phi \,\dd x'\dd t 
+ \langle \bs\res_s^h,\bs\psi\rangle\,,
\end{align}
where $\bs\res_b^\eps$ denotes the boundary residual term given by
\begin{align*}
\bs\res_b^\eps = \nu\left(\eps\pa_3v_1 -
\eps^3\int_{-1}^0(\pa_{11}v_1 + \pa_{12}v_2), \eps\pa_3v_2 -
\eps^3\int_{-1}^0(\pa_{21}v_1 + \pa_{22}v_2), 
- 2\eps^2(\pa_1v_1 + \pa_2v_2)\right)\,,   
\end{align*}
$\langle \bs\res_s^h,\bs\psi\rangle$ denotes the structure residual term $\bs\res_s^h$ acting
on a test function $\bs\psi$ as
\begin{align*}
\langle \bs\res_s^h,\bs\psi\rangle &= 
-\vro_s^h\T^{-1}\int_0^T\!\!\int_{S_h}\pa_t\bs\ua^h\cdot\pa_t\bs\psi\,\dd\bx\dd t\\
&\quad + \T\int_0^T\!\!\int_{S_h}\Big(\frac{(\lambda^h)^2}{2\mu^h + \lambda^h}\diver\bs\ua^h 
\diver(\psi_1,\psi_2) + \lambda^h\diver\bs\ua^h \pa_3\psi_3\Big)\,\dd \bx\dd t\,,
\end{align*}
and coefficients $\mu^h = {\mu}h^{-\kappa}$, $\lambda^h = {\lambda}h^{-\kappa}$ and
$\vro_s^h = {\vro}_sh^{-\kappa}$ are according to the scaling ansatz (S2).
Defining the fluid error $\bs e_f^\eps:= \bv^\eps - \bs\va^\eps$ 
and the structure error $\bs e_s^h:= \bu^h - \bs\ua^h$, and subtracting (\ref{eq:weak_fsi_app}) 
from the weak form of the original weak form, we find the variational 
equation for the errors:
\begin{align}
-\vro_f\int_0^T\!\! \int_{\Omega_\eps}\bs e^\eps_f\cdot\pa_t\bs\phi\,\dd \bx\dd t \nonumber
+ 2\nu\T\int_0^T\!\!\int_{\Omega_\eps}\sym\nabla\bs e_f^\eps:\sym\nabla\bs\phi\,\dd \bx\dd t\\
-\vro_s^h\T^{-1}\int_0^T\!\!\int_{S_h}\pa_t\bs e_s^h\cdot\pa_t\bs\psi\,\dd\bx\dd t 
+ \T\int_0^T\!\!\int_{S_h}\left(2\mu^h\sym\nabla\bs e_s^h :\sym\nabla\bs\psi
+ \lambda^h\diver\bs e_s^h \diver\bs\psi\right)\,\dd \bx\dd t \label{eq:weak_fsi_error}   \\
= \T\int_0^T\!\!\int_{\Omega_\eps}f_3^\eps\phi_3 \,\dd \bx\dd t 
- \T\int_0^T\!\!\int_{\Omega_\eps}\bs\res_f^\eps\cdot\bs\phi \,\dd \bx\dd t \nonumber
- \T\int_0^T\!\!\int_{\omega}\bs\res_b^\eps\cdot\bs\phi \,\dd x'\dd t 
- \langle \bs\res_s^h,\bs\psi\rangle\,
\end{align}
for all test functions $(\bs\phi,\bs\psi)\in \W(0,T;\Omega)$.  

Next step is a careful selection of test functions in (\ref{eq:weak_fsi_error}). First we 
choose
\begin{equation}\label{def:test_psi}
\bs\psi = \T^{-1}(\pa_t e_{s,1}^\even, \pa_t e_{s,2}^\even,\pa_t e_{s,3}^\odd)\,,
\end{equation}   
where superscripts $\even$ and $\odd$ denote even and odd components of the orthogonal decomposition
of respective functions with respect to the variable $(x_3 - h/2)$.  
Observe in (\ref{def:ua_h}) 
that, up to time dependent constants, 
components of the approximate displacement $\bs\ua^h$ are respectively odd, 
odd and even with respect to $(x_3-h/2)$.
The idea of using this particular test function comes from the fact that such $\bs\psi$
annihilates large part of the structure residual term $\bs\res_s^h$ on the right hand side in (\ref{eq:weak_fsi_error})
and the rest can be controlled (cf.~\cite[Section 4.2]{BuMu19} for details).  
Concerning the fluid part, observe that approximate solutions
do not satisfy the kinematic interface condition in the horizontal components, 
i.e.~$\va^\eps_\alpha \neq \pa_t\ua^h_\alpha$ on $\omega\times(0,T)$ and therefore
$(\bs\va^\eps,\bs\ua^h)$ does not belong to the space $\V(0,T;\Omega)$. 
For the third component however, the interface condition is satisfied.
In order to match interface values of $\bs\psi$, the fluid test function $\bs\phi$
has to be accordingly corrected fluid error, i.e.~we take
\begin{equation}   
\bs\phi = \bs e_f^\eps + \bs\varphi\,,
\end{equation}
where the correction $\bs\varphi$ satisfies 
\begin{align*}
\diver\bs\varphi &= 0\quad\text{on }\ \Omega_\eps\times(0,T)\,,\\ 
\left.\varphi_\alpha\right|_{\omega\times(0,T)} &= -\T^{-1}\left.\pa_tu_\alpha^\odd\right|_{\omega\times(0,T)},\\
\left.\varphi_3\right|_{\omega\times(0,T)} &= -\T^{-1}\left.\pa_te_{s,3}^\even\right|_{\omega\times(0,T)}\,,\\
\left.\bs\varphi\right|_{\{x_3=-\eps\}\times(0,T)} &= 0\,, 
\end{align*}
and $\bs\varphi(\cdot,t)$ is $\omega$-periodic for every $t\in(0,T)$. 
This choice of $\bs\varphi$ ensures
the kinematic boundary condition $\bs\phi = \bs\psi$ a.e.~on $\omega\times(0,T)$. 
Following \cite{BuMu19} it can be proved that 
the corrector $\bs\varphi$ satisfies the uniform bound
\begin{equation*}
\|\nabla\bs\varphi\|_{L^\infty(0,T;L^2(\Omega_\eps))} \leq C\eps^{5/2}\,
\end{equation*}
with $C>0$ independent of $\bs\varphi$ and $\eps$. Moreover, careful 
estimation of other residual terms in (\ref{eq:weak_fsi_error}) provides
the basic error estimate:
for a.e.~$t\in(0,T)$ we have
\begin{align}
\frac{\vro_f}{4}\int_{\Omega_\eps}|\bs e_f^\eps(t)|^2\dd\bx + \nonumber
\frac{\nu\T}{2}\int_0^t\!\!\int_{\Omega_\eps}|\nabla\bs e_f^\eps|^2\,\dd \bx\dd s
+\frac{\vro_s^h\T^{-2}}{4}\int_{S_h}\left((\pa_t e_{s,\alpha}^\even(t))^2 
+ (\pa_t e_{s,3}^\odd(t))^2\right)\dd\bx \\ 
+ \int_{S_h}\left(\mu^h\left|\sym\nabla ( e_{s,1}^\even,  e_{s,2}^\even, e_{s,3}^\odd)(t)\right|^2
+ \frac{\lambda^h}{2}\left|\diver( e_{s,1}^\even,  e_{s,2}^\even, e_{s,3}^\odd)(t)\right|^2 
\right)\,\dd \bx \label{ineq:error_est2}\\
 \leq C\T\eps^3(h^\gamma + h^{4\gamma - 2\kappa + 4})\,.\nonumber
\end{align}

Estimate (\ref{ineq:error_est2}) is now sufficient to conclude the error estimates for the fluid part: 
velocities and the pressure, while for the structure part, the Griso decomposition of 
the structure error needs to be examined and employing another pair of test functions 
$(\bs\psi,\bs\phi) = (\T^{-1}\pa_t\bs e_s^{\odd}, \bs e_f^\eps + \bs\varphi)$ with appropriate
corrector $\bs\varphi$ will provide sufficient conditions to conclude the error estimates also for the 
structure part (cf.~\cite[Section 4.4]{BuMu19}).

$\hfill\Box$



\subsection{Nonlinear $\eps + 0$ problem}\label{sec:nonlin}

In this section we discuss a nonlinear FSI problem in which nonlinearities appear both
in equations for the fluid motion and in geometry of the fluid domain. More precisely, the 
coupling conditions are also nonlinear and the coupling is realized on the moving interface. 
Unlike in the previous 
section, the structure is here modeled by a lower-dimensional elasticity model, and 
additionally we decrease dimensionality of the original problem to two space dimensions for the 
fluid and one space dimension for the structure. The main reason for this ad hoc dimension reduction 
in the FSI problem is availability of the wellposedness results.

Let us now describe our setting. The fluid domain at time $t$ is assumed to be of the form
\begin{equation*}
\Omega_{\eta}(t)=\{(x,z):x\in \omega,\; z\in (0,\eta(t,x))\}\subset\R^2,
\end{equation*}
where $\omega=(0,1)$, and function $\eta$ describes the dynamics of 
the vertical displacement of the top boundary.  
Let us further denote the space-time cylinder
$$
\Omega_\eta(t)\times (0,T):=\bigcup_{t\in 
(0,T)}\Omega_{\eta}(t)\times\{t\}\subset\R^2\times (0,\infty),\quad T\in (0,\infty]\,,
$$
to be domain of our problem.
The FSI problem is described by the system of partial differential equations:
\begin{align}
\vro_f(\pa_t\bv+(\bv\cdot \nabla)\bv) - \diver\sigma_f(\bv,p) &= \bs f\,,\quad 
\Omega_\eta(t)\times(0,\infty)\,,\label{1.eq:stokes}\\
\diver \bv &= 0\,,\quad \Omega_\eta(t)\times(0,\infty)\,,\label{1.eq:divfree}\\
\vro_s\pa_{tt}\eta - D \partial^2_{x}\pa_t\eta + B\partial^4_{x}\eta  
&=-J^{\eta} \big (\sigma_f(\bv,p){\bf n}^{\eta}\big )(t,x,\eta(t,x))\cdot{\bf e}_z\,, 
\quad \omega\times(0,\infty)\,,\label{1.eq:elast}
\\
\bv(t,x,\eta(t,x))&=(0,\partial_t\eta(t,x))\,,\quad \omega\times(0,\infty)\,,\label{KinematicBC}
\end{align}
where $\sigma_f(\bv,p) = 2\nu \sym \nabla \bv - pI$ denotes the Cauchy stress tensor of the 
viscous fluid, $\nu$, $\varrho_f>0$ are the fluid viscosity and density, respectively, and 
$\bs f$ denotes the fluid external force.
Furthermore, $\varrho_s$ is 
the structure density, $\bs n^{\eta}$ is the unit outer normal to the deformed configuration 
$\Omega_{\eta}$, $J^{\eta}(t,x)=\sqrt{1+\partial_x\eta(t,x)^2}$ is Jacobian of the transformation 
from Eulerian to Lagrangian coordinates, and constants $D,B>0$ describe visco-elasticity and 
elasticity properties of the structure, respectively.

Equations \eqref{1.eq:stokes} and \eqref{1.eq:divfree} are standard incompressible 
Navier-Stokes equations describing the flow of the Newtonian fluid, while the structure is 
described by a linear equation of visco-elastic plate (\ref{1.eq:elast}). 
The fluid and the structure are coupled via dynamic and kinematic coupling 
conditions \eqref{1.eq:elast} and \eqref{KinematicBC} representing the balance of forces 
in ${\bf e}_z$ direction and continuity of the velocity, respectively.
Additional simplifying assumption is that the structure moves only in the vertical direction.
This is not fully justified from the physical grounds, but it is reasonable in the view of
results of the previous section, where it is shown that, up to time-dependent translations,
the bending regime is the dominant one and the displacement in the horizontal direction is 
of the lower order.
For more details about physical background of system \eqref{1.eq:stokes}-\eqref{KinematicBC} and 
corresponding lower-dimensional elasticity models we refer to \cite{BorSun,Cia97} and 
reference therein. 
The bottom boundary is rigid and we prescribe the standard no-slip boundary condition for the 
fluid velocity:
$\bv(t,x',0)=0$ for all $(x',t)\in\omega\times(0,\infty)$.
On the lateral boundaries we prescribe the periodic boundary conditions in the horizontal direction,
which is taken for technical simplicity and because of availability of the 
global existence results \cite{grandmont2016existence}. 
In such a case the flow is driven by the right-hand side $\bs f$.
Finally, for simplicity of exposition, we impose some trivial initial conditions:
$\bv(0,.)=0,\; \eta(0,.)=\eta_0,$ and $\partial_t\eta(0,.)=0$,

Since our aim is to derive the reduced model in the regime of relatively thin domain, we assume 
that the initial thickness of the domain is $O(\eps)$, i.e.~$\eta^{\eps}_0(x)=\eps\eta_0(x)$ for 
$x\in (0,L)$. 
Moreover, like in the previous section, we assume that 
$\|\bs f\|_{L^\infty(0,\infty;L^{\infty}(\Omega_\eta(t);\R^2))} \leq C$, which is satisfied by
physically relevant volume forces. 
Testing formally equations (\ref{1.eq:stokes}) and (\ref{1.eq:elast}) with classical solutions
$\bv$ and $\pa_t\eta$,
respectively, and integrating by parts yields the basic energy inequality:
for every $t\in (0,T)$
\begin{align}
\frac{\varrho_f}{2}\|\bv(t)\|^2_{L^2(\Omega_{\eta}(t))}\nonumber
+2\nu\int_0^t\!\!\int_{\Omega_{\eta}(s)}|\sym\nabla\bv|^2\dd \bx\dd s\\
+\frac{\varrho_s}{2}\|\partial_t\eta(t)\|^2_{L^2(\omega)} \label{EI}
+D\int_0^t\|\partial_t\partial_{x}\eta\|^2_{L^2(\omega)}\dd s
+\frac{B}{2}\|\partial^2_x\eta(t)\|^2_{L^2(\omega)}\\
\leq \frac{B}{2}\|\partial_x^2 \eta_0\|^2_{L^2(\omega)}\nonumber
+\int_0^t\!\!\int_{\Omega_{\eta}(s)}\bs f\cdot\bv\,\dd\bx\dd s\,. 
\end{align}

Let us now discuss weak solutions. First we introduce appropriate solution spaces. 
The fluid solution space will depend on the displacement $\eta$. If we denote
\begin{equation*}
{V}_F(t)=\left\{\bv\in H^1(\Omega_{\eta}(t))\ :\ \diver\bv=0,\;
 \bv|_{z=0}=0,\ \bv \text{ is } \omega-\text{periodic in }x_1\right\},
\end{equation*}
then the above energy estimate suggests that appropriate fluid solution space
is 
\begin{equation*}
{\mathcal V}_F(0,\infty;\Omega_{\eta}(t))
= L^{\infty}(0,\infty;L^2(\Omega_{\eta}(t)))\cap L^2(0,\infty;{V}_F(t))
\end{equation*}
while for the structure, again based on the energy estimate, we choose the solution space to be 
\begin{equation*}
{\mathcal V}_S(0,\infty;\omega)
= W^{1,\infty}(0,\infty;L^2(\omega))\cap 
L^{\infty}(0,\infty;H^2_{\rm per}({\omega}))\cap H^1(0,\infty;H^1(\omega))\,.
\end{equation*}

\begin{definition}
We call $(\bv,\eta)\in{\mathcal V}_F(0,\infty;\Omega_{\eta}(t))\times {\mathcal V}_S(0,\infty;\omega)$ a 
{\em weak solution} of the FSI problem (\ref{1.eq:stokes})-(\ref{KinematicBC}) 
if for every $(\vphi,\psi)\in C^1_c([0,\infty);\mathcal{V}_F(t)\times H^2_{\rm per}({\omega}))$ 
satisfying $\vphi(t,x,\eta(t,x'))=\partial_t\psi(t,x){\bf e}_3$ it holds
\begin{align}
-\varrho_f\int_0^\infty\!\!\int_{\Omega_\eta(t)}\left(\bv\cdot\partial_t\vphi
+(\bv\cdot\nabla)\bv\cdot\vphi\right)\dd \bx\dd t
+2\nu\int_0^\infty\!\!\int_{\Omega_\eta(t)}\sym\nabla\bv:\sym\nabla\vphi\dd\bx\dd t
\nonumber
\\
-\varrho_s\int_0^\infty\!\!\int_{\omega}\partial_t\eta\partial_t\psi\dd x_1\dd t
+D\int_0^\infty\!\!\int_{\omega}\partial_{x}\partial_t\eta\cdot\nabla\psi \dd x_1\dd t
+B\int_0^\infty\!\!\int_{\omega}\partial^2_{x}\eta\partial^2_{x}\psi \dd x_1\dd t
\label{WeakFormulation}
\\
=\varrho_s\int_{\omega}\eta_0\psi(0)\dd x_1 
+ \int_0^{\infty}\!\!\int_{\Omega_\eta(t)}\bs f\cdot\vphi\,\dd\bx\dd t
\nonumber
\end{align}
and \eqref{KinematicBC} is satisfied in the sense of traces. Moreover, the energy inequality 
(\ref{EI}) is satisfied.
\end{definition}
The existence of weak solutions is by now well-established in the literature, see e.g. \cite{CDEM,BorSun}. 
However, the existence results are not global and state that weak solutions exist as long as there 
is no contact between the elastic and rigid boundary, i.e.~in our notation as long as $\eta>0$. 
Even though there are results that contact will not occur in the case when the structure is 
rigid \cite{HilTak}, to the best of our 
knowledge there are no global in time existence results for weak solution to problem 
\eqref{1.eq:stokes}-\eqref{KinematicBC}. Since we are interested in the long time behavior, we rely 
on the following recent result on the existence of global-in-time strong solutions.

\begin{theorem}[\cite{grandmont2016existence}]\label{ExistenceGlobal}
For every $\eps>0$ there exists global-in-time strong solution to problem 
\eqref{1.eq:stokes}-\eqref{KinematicBC} with initial conditions 
$\bv_0=0$, $\eta_0=\eps\hat{\eta}_0$, $\eta_1=0$ and the right hand 
side $\bs f\in L^2(0,T;L^2_{\rm per}(\R^2))$. 
Solutions satisfy $\eta^\eps>0$ and additional regularity properties which we briefly write
\begin{equation}\label{StrongSolReg}
\eta^\eps\in H^2_tL^2_x\cap L^2_tH^4_x\,,\quad
\bv^\eps\in L^2_tH^2_x\quad\text{and}\quad
p^\eps\in L^2_tH^1_x\,.
\end{equation}
\end{theorem}

Motivated by the results from previous sections we assume the following scaling 
ansatz in order to be in the thin-film regime:

\begin{enumerate}
  \item[(S1)] $B = \hat{B}\eps^{-1}$, $D=\hat{D}\eps^{-2}$ and $\vro_s = \hat{\vro}_s\eps^{-1}$ for some
 $\hat{B},\;\hat{D},\;\hat{\vro}_s>0$ independent of $\eps$;
  \item[(S2)] $\T = \eps^{-2}$.
\end{enumerate}
To the best of our knowledge rigorous derivation of equation \eqref{2:stf} as the reduced model for
the FSI system \eqref{1.eq:stokes}-\eqref{KinematicBC} is still missing from the literature. 
Here we present the main steps of the derivation without proofs. The details of the proofs will be 
included in a forthcoming work \cite{BuMu20}. The main steps are analogous to the main steps in 
linear case, but technically much more involved. The main difficulties are consequence of the 
fact that $\eps$-problems are moving boundary problems, i.e. we have to deal with the geometrical 
nonlinearity at every step of the derivation.

\textit{Step 1: Uniform energy estimate.}
The first step is to quantify the energy estimate (\ref{EI}) in terms of the small parameter $\eps$. 
Unlike in linear case, the domain depends on solution and therefore one has to carefully track 
dependence of functional inequalities (e.g.~the Poincar\'e inequality) on the solution itself. 
By using the scaling ansatz (S1) we arrive to the following uniform (in small parameter $\eps$) 
energy inequality: for a.e.~$t\in (0,T)$
\begin{equation}\label{EIEpsFinal}
\frac{\varrho_f}{2}\|\bv^{\eps}(t)\|^2_{L^2(\Omega_{\eta}(t)}
+\nu\int_0^t\int_{\Omega_{\eta}(t)}|\nabla\bv^{\eps}|^2
+\frac{D}{\eps^2}\int_0^t\|\partial_t\partial_{x}\eta^{\eps}\|^2_{L^2(\omega)}
+\frac{B}{\eps}\|\partial^2_x\eta^{\eps}\|^2_{L^2(\omega)}
\leq  C t\eps^3\,.
\end{equation}
By taking into account scaling ansatz (S2) and using the Sobolev embedding we get an 
$L^{\infty}$-estimate for the displacement: 
\begin{equation}\label{DisplacementEstimate}
\|\eta^{\eps}\|_{L^{\infty}(0,T;L^{\infty}(0,L))}\leq C\eps\,.
\end{equation}

\textit{Step 2: Positivity of the limit displacement.} Let us denote by $\eta$ the weak limit of 
$\eta^{\eps}/\eps$. From Theorem \ref{ExistenceGlobal} it is immediate that $\eta\geq 0$. However, 
in order to perform our analysis we need to prove the strict positivity $\eta>0$. 
This can be done by adapting estimates used in the proof of 
Theorem \ref{ExistenceGlobal} \cite{grandmont2016existence} to our case and combining them 
with the scaling ansatz.

\textit{Step 3: ALE formulation.}
In order to identify the limit model we need to reformulate weak formulation \eqref{WeakFormulation} on the fixed reference domain. The main difference in comparison to the linear case is that now the change of variable depends on the solution itself. In numerical computation this formulation is usually called Arbitrary Lagrangian-Eulerian (ALE) formulation. We use the following explicit form of the change of variables:
\begin{equation}\label{VarChange}
\left (
\begin{array}{c}
\hat{t}\\
y_1\\
y_2
\end{array}
\right )=
\left (
\begin{array}{c}
\eps^2 t\\
x_1\\
\frac{x_2}{\eta^{\eps}(t,x_1)}
\end{array}
\right ).
\end{equation}

\textit{Step 4: Identifying the limit model.}
In the last step we pass to the limit as $\eps\downarrow0$ and obtain the limit model \eqref{2:stf}. 
The limiting procedure follows heuristic described in the introduction and is similar as in the 
linear case. The main difference is that due to the nonlinearities in the system, we need to prove 
the strong convergence properties of sequence $(\bv^{\eps},\eta^{\eps})$ in order to pass to the limit. 
For this we need two main ingredients: Aubin-Lions lemma for strong convergence of the displacement 
sequence $\eta^{\eps}$ and the scaling ansatz for the convergence of the fluid convective term.

\subsection*{Acknowledgement}
This work has been supported in part by the Croatian Science
Foundation under projects 7249 (MANDphy) and 3706 (FSIApp).

\end{document}